\documentclass[onefignum,onetabnum]{siam_latex_template/siamart190516}
\usepackage[utf8]{inputenc}
\usepackage{textcomp}

\usepackage{graphicx}
\usepackage{multicol}
\usepackage{amsmath}
\usepackage[version=4]{mhchem}
\usepackage{siunitx}
\usepackage{longtable,tabularx}

\usepackage{lipsum}
\usepackage{tikz}
\usepackage{subfig}
\usepackage{dsfont}
\usepackage{amsfonts}

\usetikzlibrary{shapes,arrows}

\title{Resource-Constrained Optimal Experimental Design}

\author{Anthony M. DeGennaro\thanks{Associate Computational Scientist, Computational Science Initiative, Brookhaven National Laboratory (\email{adegennaro@bnl.gov}).}
\and Francis J. Alexander\thanks{Deputy Director, Computational Science Initiative, Brookhaven National Laboratory.}}

\begin{document}

\tikzstyle{block} = [draw, fill=blue!20, rectangle, 
    minimum height=3em, minimum width=6em]
\tikzstyle{sum} = [draw, fill=blue!20, circle, node distance=1cm]
\tikzstyle{input} = [coordinate]
\tikzstyle{output} = [coordinate]
\tikzstyle{pinstyle} = [pin edge={to-,thin,black}]

\maketitle

\begin{abstract}
  
  The goal of this paper is to make Optimal Experimental Design (OED)
  computationally feasible for problems involving significant
  computational expense. We focus exclusively on the Mean Objective
  Cost of Uncertainty (MOCU), which is a specific methodology for OED,
  and we propose extensions to MOCU that leverage surrogates and
  adaptive sampling. We focus on reducing the computational expense
  associated with evaluating a large set of control policies across a
  large set of uncertain variables. We propose reducing the
  computational expense of MOCU by approximating intermediate
  calculations associated with each parameter/control pair with a
  surrogate. This surrogate is constructed from sparse sampling and
  (possibly) refined adaptively through a combination of sensitivity
  estimation and probabilistic knowledge gained directly from the
  experimental measurements prescribed from MOCU. We demonstrate our
  methods on example problems and compare performance relative to
  surrogate-approximated MOCU with no adaptive sampling and to full
  MOCU. We find evidence that adaptive sampling does improve
  performance, but the decision on whether to use
  surrogate-approximated MOCU versus full MOCU will depend on the
  relative expense of computation versus experimentation. If
  computation is more expensive than experimentation, then one should
  consider using our approach.
  
\end{abstract}




\section{Introduction}

It is frequently the case that engineers and scientists
must make decisions under uncertainty and with only partial
information available. This is especially so when designing complex
systems. For one, the space of all designs or controls that could
possibly be selected (the design space) usually has an enormously
large parametric dimension. For another, the process used to evaluate
any single design choice is often very costly, measured either by
required time or resources. Lastly, complex engineering designs are
usually situated in complex physics. Here, operating conditions,
initial/boundary conditions, and various parametric models must be
known and specified, but it is never possible to do this with complete
accuracy or precision. These uncertainties can propagate through the
physics and strongly affect the output quantities-of-interest (QOIs),
even if their magnitude is relatively small (as is the case, for
example, in the state-space trajectories of chaotic systems). In all
cases, the expense associated with the needed computation forces
researchers to operate with knowledge of only a subset of the full
design space. The methods we use to explore and to optimize over these
design spaces must take these constraints into consideration; if they
do not, they are not helpful for ``real-world'' problems.

The need for efficient, calculated exploration of parametric spaces
and judicious allocation of resources in engineering design problems
has given rise to the field of Optimal Experimental Design (OED). OED
can be conceptualized as a combination of optimization under
uncertainty (OUU) with Bayesian calibration to experimental data. The
goal is to design a system that is optimal with respect to a given
goal, on average across some uncertain parameters. Experiments may be
conducted that could reveal with greater statistical
accuracy/precision what the true values of the uncertain parameters
are, but these experiments are costly to perform. Therefore, the
second objective of OED is to select a series of experiments that are
the most informative with respect to the engineering goal. The result
is a scheme in which data-driven calibration and optimization
co-depend on each other. Intelligently selected experiments
iteratively reduce uncertainty, allowing one to optimize over the
design space with greater certainty and accuracy. At the same time,
those optimization results are used to suggest which experiments would
be most useful to conduct, as measured by the given engineering
goal. The Mean Objective Cost of Uncertainty
(MOCU)~\cite{boluki2018experimental, yoon2013quantifying,
  imani2018sequential, dehghannasiri2014optimal,
  dehghannasiri2015efficient, yoon2020quantifying} is one specific
algorithm that implements OED; other alternatives include
entropy-based exploration strategies, active learning, and Knowledge
Gradient~\cite{frazier2009knowledge, frazier2008knowledge,
  ryzhov2012knowledge, chen2013optimistic}.

We propose in this paper an extension to MOCU that gives users
a principled way of dealing with constraints on computational resources and sparse
data. The goal of this is to make MOCU better suited to the demands of
design problems that are computationally intensive. Our approach is to
first build a surrogate model that approximates the cost of any
particular design choice. This model is coarse, in the sense that it
is trained on a sparse data set drawn from the full design space. This
model is then updated with adaptively selected sample points as more
information is gathered about the design space through experimental
measurements. In other words, we are proposing an approximate OED
scheme, where the accuracy of the surrogate used to drive MOCU
co-evolves with, and depends on, the experimental knowledge obtained
from doing MOCU. In MOCU, experiments are chosen both to reduce
uncertainty about the design space and to optimize some QOI over that
space. Our insight is that, when using an imperfect surrogate model
for how the QOI depends on the design space, those experiments also
carry information that should be used to optimally refine the
surrogate. To our knowledge, this contribution represents one of the
first efforts at making MOCU/OED possible for large and complex
systems where computational resources are constrained.

This paper is organized as follows. We begin with a brief review of
OED/MOCU, as well as some of the tools we will be using for adaptive
sampling and surrogate construction. We then proceed to a description
of our approach for sparse, adaptive MOCU. We close with several
examples that explicitly compare our methods to standard MOCU and to
sparse MOCU with no adaptive sampling. We conclude that adaptive
sampling does improve MOCU performance relative to static surrogates,
but that the decision on whether to use surrogate-approximated MOCU
versus full MOCU will depend on the relative expense of computation
versus experimentation. If computation is expensive relative to
experimentation, then one should consider using our approximation
methods. This might be the case, for example, when multiscale physics
codes are involved, and must be evaluated over many parametric
settings, operating conditions, and initial/boundary conditions. If,
instead, experimentation is the expensive step, then one should
consider traditional full MOCU. This might occur, for example, if
doing a single experiment involves fabricating and testing a new
exotic material in a laboratory setting.

\section{Background}

In order to make this paper as self-contained as possible, we provide
in this section a brief review of the literature for important
algorithms we will make use of. We begin with an overview of the MOCU
algorithm, which is central to this research. We then give a short
overview of Gaussian Process (GP) regression, which is our tool of
choice for surrogate construction.

\subsection{The MOCU Algorithm}

The MOCU algorithm is an approach to OED that ties together Bayesian
calibration with OUU. The overall goal of MOCU is to minimize a cost
function of two variables, $J(\theta,\psi)$. Here, $\theta$ is an
element of a discrete set, $\Theta = \lbrace \theta_1 , \dots ,
\theta_{n_{\theta}} \rbrace$, called the uncertainty class. For
example, $\Theta$ could consist of $n_{\theta}$ independent draws from
some underlying continuous probability distribution. $\psi$ is also an
element of a discrete set, $\Psi = \lbrace \psi_1 , \dots ,
\psi_{n_{\psi}} \rbrace$, called the action set. The distinction
between $\theta$ and $\psi$ is control: $\psi$ is a variable which we
can directly control, while $\theta$ is an uncertain parameter, whose
possible values are given by the set $\Theta$.  We also assume we have
a prior belief about the most probable elements of $\Theta$, which we
have in the form of a probability mass function $\rho(\theta)$. The
cost function $J(\theta,\psi)$ quantifies loss associated with a
design objective. To give a concrete example from aerodynamics:
$J(\theta,\psi)$ may be the inverse of the lift-to-drag ratio for a
given airfoil, which we wish to minimize. $\Psi$ would consist of a
set of airfoil shapes, and $\Theta$ could consist of different
parameter values for the turbulence closure model used by the
numerical flow solver that computes the flow over the airfoil.

If we knew the ground truth value of $\theta$, which we will call
$\theta_{\text{true}}$, then we could simply solve a one-parameter
optimization problem:

\begin{equation}
\label{eq:psi_opt}
\psi_{\theta_{\text{true}}} = \text{argmin}_{\Psi} \; J(\theta=\theta_{\text{true}},\psi) \;\;\; .
\end{equation}
Of course, we do not know $\theta_{\text{true}}$, and so the best we
can do is select that control policy $\psi_{\rho(\theta)}$ which
minimizes $J(\theta,\psi)$ over the distribution of $\theta$:

\begin{equation}
\label{eq:psi_mocu}
\psi_{\rho(\theta)} = \text{argmin}_{\Psi} \; \mathbb{E}_{\rho(\theta)} [ J(\theta,\psi) ] \;\;\; .
\end{equation}

Up to this point, this problem formulation is no different than
standard OUU. However, we now further assume that we have the ability
to conduct experiments. The measurements that we make from an
experiment carry information about the what elements of $\Theta$ are
statistically most likely to be true. There are a number of possible
experiments that we could conduct, and we denote that set as $X =
\lbrace x_1 , \dots , x_{n_X} \rbrace$. The set of possible outcomes
for any particular experiment $x_i$ will be denoted $Y = \lbrace y_1 ,
\dots , y_{n_Y} \rbrace$.

This addition of experiments makes the problem different from standard
OUU because the posterior distribution $\rho(\theta | x,y)$ will
depend on what specific experiment $x$ is selected, and what the outcome
$y$ of that experiment is (once conducted). Now, we should optimize
the cost function over the experiment-conditioned posterior (rather
than the prior as in Eq.~\ref{eq:psi_mocu}):

\begin{equation}
\label{eq:psi_mocu_posterior}
\psi_{\rho(\theta|x,y)} = \text{argmin}_{\Psi} \; \mathbb{E}_{\rho(\theta|x,y)} [ J(\theta,\psi) ] \;\;\; .
\end{equation}

Eq.~\ref{eq:psi_mocu_posterior} yields an uncertainty-robust, optimal
control policy for each possible experiment and its possible
outcomes. While $\psi_{\rho(\theta|x,y)}$ is the best strategy ``on
average'', it is not likely to be the optimal strategy for any
particular element of $\Theta$ (such as, in particular, the actual
value $\theta_{\text{true}}$). We denote the $\theta$-specific optimal
strategy as $\psi_{\theta}$. To determine which experiment should be
conducted, we need to consider the cost of using
$\psi_{\rho(\theta|x,y)}$ rather than $\psi_{\theta}$, averaged over
the experiment-conditioned posterior:

\begin{equation}
\label{eq:psi_mocu_posterior_2}
M_{\Psi}(\theta|x,y) \equiv \mathbb{E}_{\rho(\theta|x,y)} [ J(\theta,\psi_{\rho(\theta|x,y)}) - J(\theta,\psi_{\theta}) ] \;\;\; .
\end{equation}

This quantity is known as the mean objective cost of uncertainty. We
wish to select that experiment which minimizes this quantity, averaged
over all potential experimental outcomes:

\begin{equation}
\label{eq:opt_exp}
x^* = \text{argmin}_{X} \; \mathbb{E}_{y} [ M_{\Psi}(\theta|x,y) ] \;\;\; .
\end{equation}

In summary, Eq.~\ref{eq:opt_exp} prescribes the experiment that should
be conducted, and Eq.~\ref{eq:psi_mocu_posterior} gives the
uncertainty-robust, optimal control strategy. Importantly, we note
that in order to solve these equations, one must evaluate the cost
function $J(\theta,\psi)$ over all $(\theta,\psi) \in \Theta \times
\Psi$. Accordingly, the focus of this paper will be on performing MOCU
when $J(\theta,\psi)$ is approximated by a surrogate.

\subsection{Gaussian Process Regression}

Function approximation is a field that attempts to replace an
input-output mapping with a surrogate that can be queried at any input
with reduced computational expense. Further, in situations where
the input space is of high dimension and/or the evaluation of any
particular input is computationally expensive, it should be possible
to construct this surrogate with data that is relatively sparse.

Many techniques exist for this purpose. Classical methods revolve
around linear regression and spectral methods, while modern tools
include Gaussian Process Regression, Polynomial Chaos
Expansions~\cite{ghanem1990polynomial, xiu2002wiener}, and machine
learning. We focus exclusively in this paper on GP
regression~\cite{williams2006gaussian}. GP regression can be
conceptualized as a statistical technique wherein the goal is to learn
a distribution of functions that best fit the training data. One
begins by defining a prior distribution over candidate functions. This
distribution is over the parameters of a Gaussian Process, and it
defines probabilistic ranges for the mean and covariance of any
function drawn from it. Training data is then collected, and this is
used to generate a posterior distribution over the space of candidate
functions via Bayes' theorem. This means that the result of GP
regression is a distribution of surrogates, rather than just a single
point-estimate. This makes it possible to compute error bounds related
to exploratory ignorance and the functional topography for any given
input. We note that high parameter space dimension poses problems for
GP regression. However, we will only be using GP regression to
approximate a function of two variables in this work.

\section{Approach}

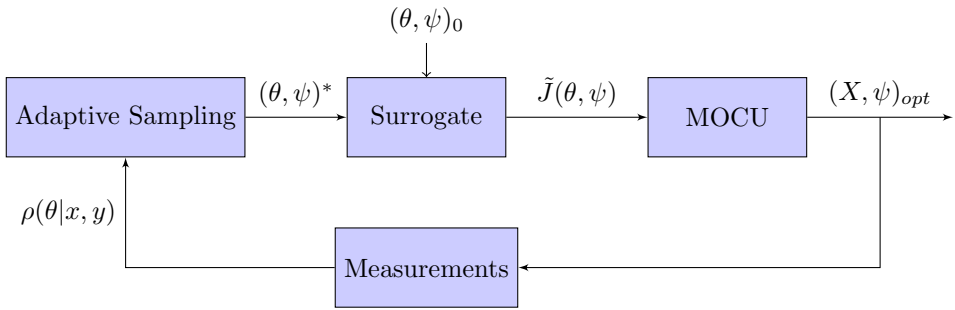
\begin{figure}[h!]
\centering
\begin{tikzpicture}[auto, node distance=2cm,>=latex']
    \node [block, 
      node distance=3cm] (samplepts) {Adaptive Sampling};
    \node [block, right of=samplepts, pin={[pinstyle]above:$(\theta,\psi)_0$},
      node distance=4cm] (surrogate) {Surrogate};
    \node [block, right of=surrogate,
      node distance=4cm] (mocu) {MOCU};
    
    \draw [->] (samplepts) -- node[name=u] {$(\theta,\psi)^*$} (surrogate);
    \node [output, right of=mocu, node distance=3cm] (output) {};
    \node [block, below of=surrogate] (measurements) {Measurements};

    \draw [->] (surrogate) -- node {$\tilde{J}(\theta,\psi)$}(mocu);
    \draw [->] (mocu) -- node [name=y] {$(X,\psi)_{opt}$}(output);
    \draw [->] (y) |- (measurements);
    \draw [->] (measurements) -| node[pos=0.99] {} 
        node [near end] {$\rho(\theta | x,y)$} (samplepts);
\end{tikzpicture}

\caption{Schematic for sparse, adaptive MOCU.}
\label{fig:mocu_tikz}
\end{figure}

A schematic depicting our approach is shown in
Fig.~\ref{fig:mocu_tikz}. As has been noted, the most computationally
intensive stage of MOCU involves computing the cost $J(\theta,\psi)$
for all $(\theta,\psi) \in \Theta \times \Psi$. To reduce the
computational expense of this, we use a surrogate model to approximate
the cost matrix $\mathcal{J} \in \mathbb{R}^{n_{\theta} \times
  n_{\psi}}$. This surrogate is constructed initially using a sparsely
sampled data set $\mathcal{P} = (p_1 , \dots , p_s) \; , \; p_i =
(\theta,\psi)_i \in \Theta \times \Psi$. This initial set of $s$
sample points -- denoted $(\theta,\psi)_0$ in Fig.~\ref{fig:mocu_tikz}
-- are drawn uniformly in $\psi$ and from the prior distribution
$\rho(\theta)$ in $\theta$. The surrogate is constructed using
Gaussian Process (GP) Regression, which is implemented in the Scipy
Sklearn library. We use a Matern kernel, and optimize the kernel
hyperparameters using gradient ascent on the log-marginal-likelihood
function. We use the mean values of the resulting GP predictions at
all locations of the cost matrix, resulting in the approximation
$\widetilde{\mathcal{J}}$.

The approximate cost matrix is fed into the MOCU algorithm, which
determines the optimal experiment, optimal policy, and a set of
conditional posteriors on $\theta$ (one for each possible outcome of
the experiment). The experiment $x$ is then conducted, and the result
$y$ of this experiment informs the conditional posterior $\rho(\theta
| x,y)$. In our algorithm, this posterior is fed back into the MOCU
algorithm along with the approximate cost matrix (just as in standard
MOCU) until the variance of $\rho(\theta | x,y)$ has fallen below a
certain preset fraction of its initial value. We note that this is a
parameter that must be selected by the user, and in general it will
depend on the details of the problem -- particularly on the likelihood
function, which plays a large role in determining the convergence of
the posterior.

Once the posterior has converged sufficiently as described, we compute
the inner 68th percentile range of the posterior, and locate the
subset of training points $\mathcal{P}_{68} \in \mathcal{P}$ within
that range. The usage of the inner 68th percentile is a heuristic
based on the amount of probability mass contained within $\pm 1\sigma$
for a normal distribution. We then test the sensitivity of the GP
model predictions to each of the points in $\mathcal{P}_{68}$. We do
this through ``leave-one-out'' validation. By this, we simply mean
that the GP sensitivity to training point $p_i$ is computed by
omitting $p_i$ from $\mathcal{P}$, re-computing a new GP surrogate
from that down-sampled training set, and computing the
$\mathcal{L}_2$-norm of the error between the new cost matrix
predictions and the original. Once we have located that training point
with the highest sensitivity, we add new training points in the
vicinity of it, sampled at random from a Gaussian distribution
centered at that point (these are denoted $(\theta,\psi)^*$ in
Fig.~\ref{fig:mocu_tikz}). We then compute a new surrogate with this
augmented training set, and continue with MOCU as before.

\section{Examples and Results}

We now proceed to apply surrogate-approximated MOCU to some example
problems, in order to assess its performance in locating the optimal
policy $\psi_{\theta_{\text{true}}} = \text{argmin}_{\Psi} \;
J(\theta_{\text{true}} , \psi)$. To this end, we compute results for
each example problem once using an adaptively refined surrogate, a
second time using a static/non-refined surrogate, and a third time
using full MOCU with no surrogates. We then compare how successful
each method was in correctly locating the optimal policy. We also
address other related topics, such as the rate of convergence
of the uncertainty class distribution and how that is affected by the
approximation scheme, and the qualitative properties of the method we
use for computing sensitivity and adaptive refinement. In our
examples, we show how our approximate MOCU methods may be used in a
multifidelity setting, and how they may be applied to the design and
control of physical systems.

\subsection{Single Model Cost Function}

Our goal here is to demonstrate our method on a fabricated cost function:

\begin{equation}
\begin{aligned}
J(\theta,\psi) &= J_1(\theta,\psi) \left(1 - J_2(\psi) J_3(\theta) \right) \\
J_1(\theta,\psi) &= 2 - \text{exp}\left[ -\frac{1}{2} \left( \psi - \frac{\theta^2}{2 n_{\theta}^2} \right)^2 / \left(\frac{n_{\theta}+n_{\psi}}{8}  \right)^2\right] \\
J_2(\psi) &= \text{exp}\left[ -\frac{1}{2} \left(\psi - \frac{3}{4}n_{\psi}\right)^2 / \left(\frac{n_{\psi}}{16}\right)^2 \right] \\
J_3(\theta) &= \text{exp}\left[ -\frac{1}{2} \left(\theta - \frac{1}{4}n_{\theta}\right)^2 / \left(\frac{n_{\theta}}{8}\right)^2 \right] \\
\Theta &= \lbrace 1 , \dots , n_{\theta}  \rbrace \; , \; \Psi = \lbrace 1 , \dots , n_{\psi}  \rbrace \; , \; \theta_{\text{true}} = \frac{1}{4} n_{\theta}
\end{aligned}
\end{equation}

\begin{figure}[h!]
  \begin{minipage}{0.32\textwidth}
    \centering 
    \subfloat[Cost function]{ \includegraphics[trim=4cm 0cm 5cm 0cm, clip, width=0.90\linewidth]{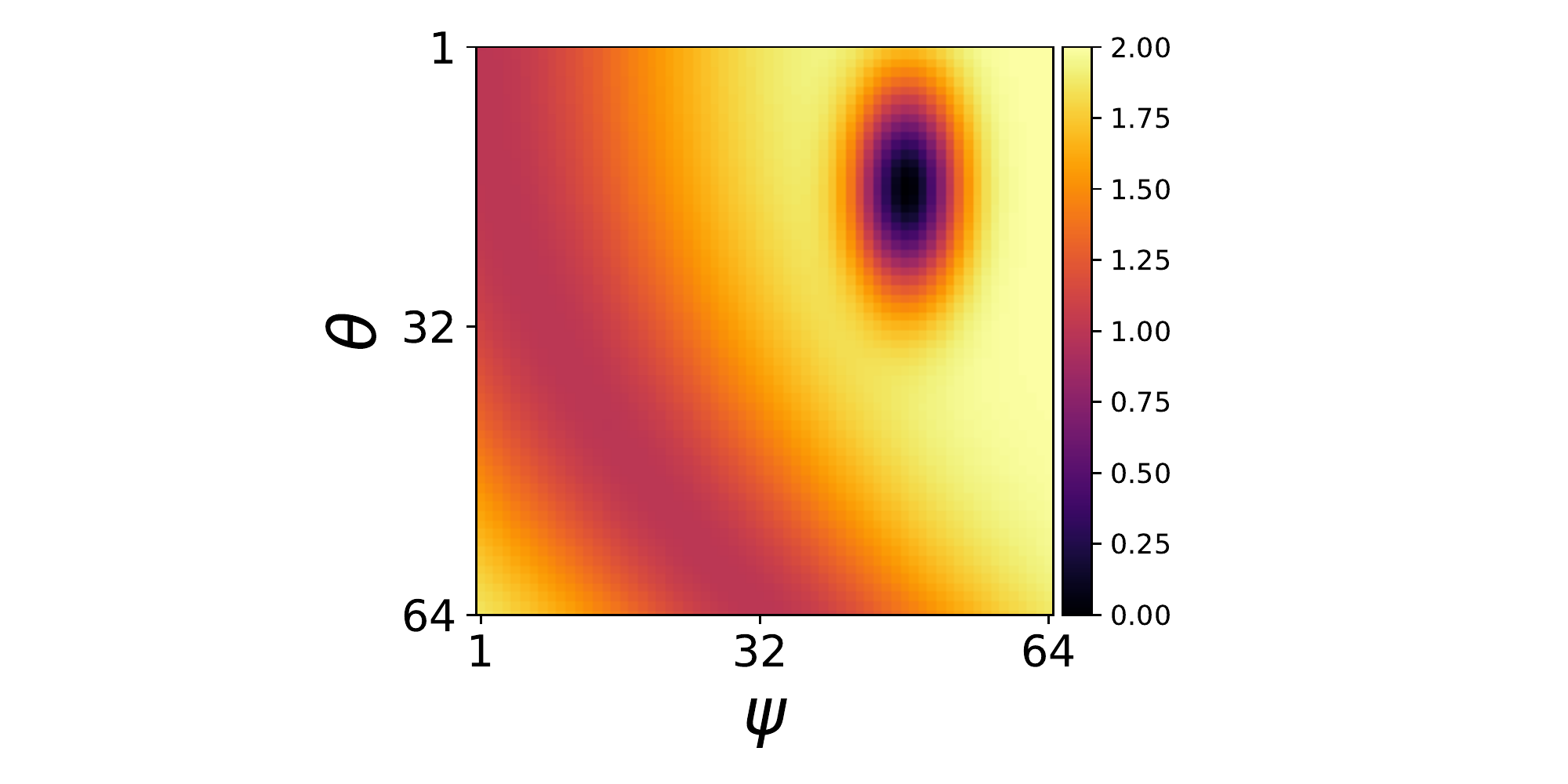} \label{fig:costfunction_singlemodel_topo} }
  \end{minipage}
  \begin{minipage}{0.32\textwidth}
    \subfloat[MC gradient]{ \includegraphics[trim=4cm 0cm 5cm 0cm, clip, width=0.90\linewidth]{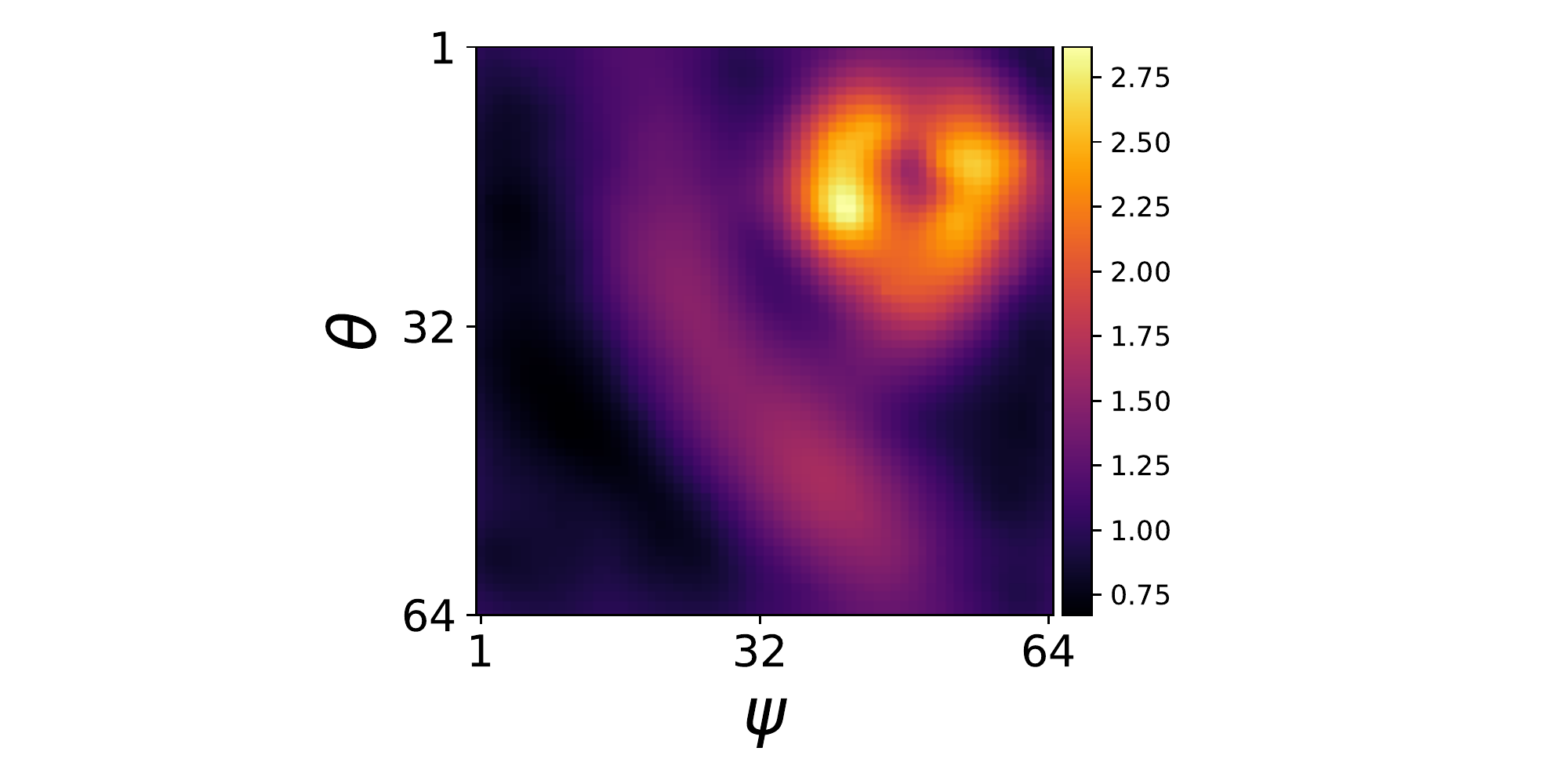} \label{fig:costfunction_singlemodel_DJ} }
  \end{minipage}
  \begin{minipage}{0.32\textwidth}
    \subfloat[MC Sensitivity]{ \includegraphics[trim=4cm 0cm 5cm 0cm, clip, width=0.90\linewidth]{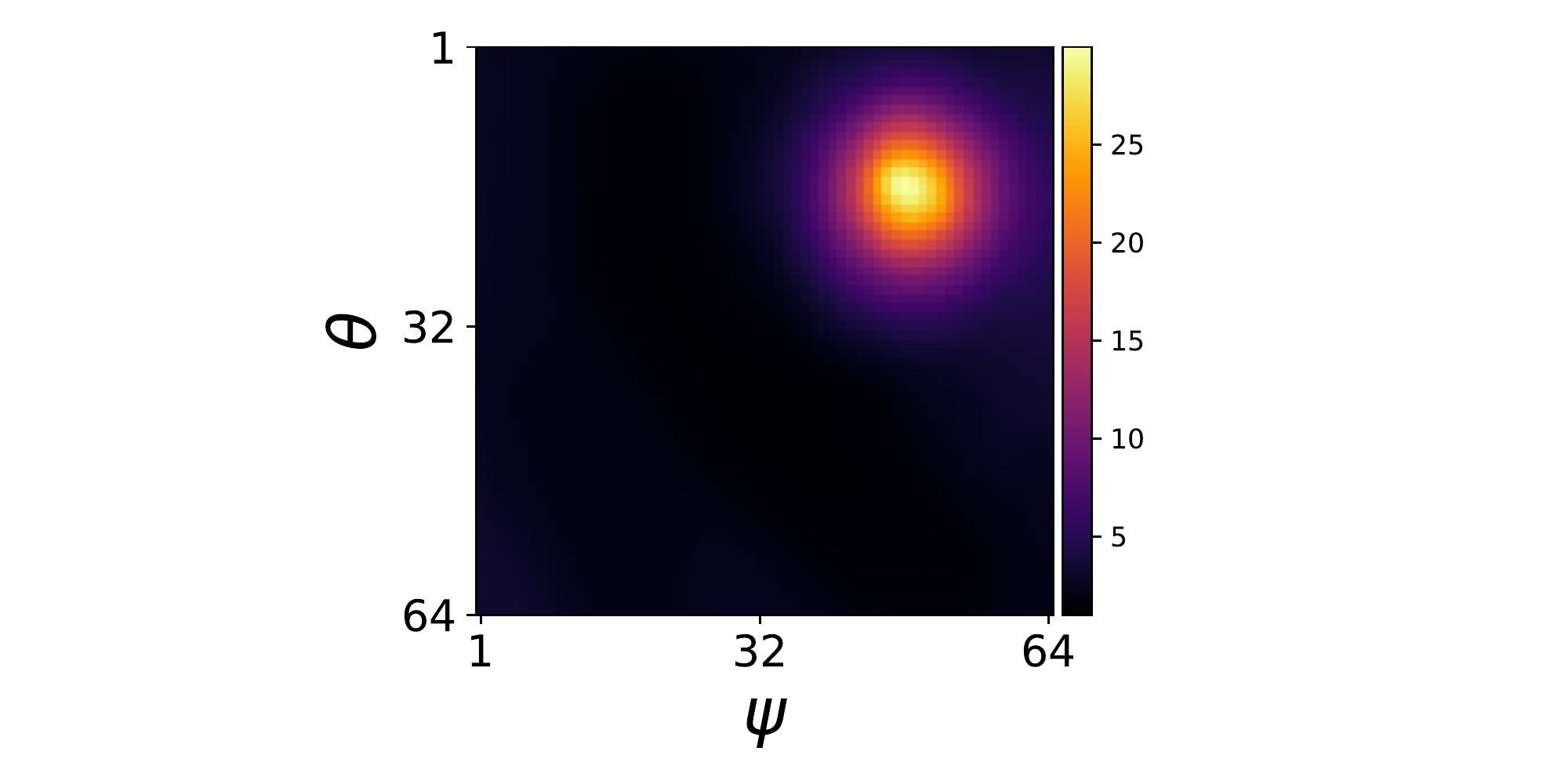} \label{fig:costfunction_singlemodel_sensitivity} }
  \end{minipage}
\caption{(a) Topography of cost function, and GP estimates of its (b) gradient magnitude and (c) local sensitivity. (b) and (c) are averages over 64 independent realizations, where one realization corresponds to a random draw of 48 initial training data points for the construction of the GP surrogate of (a).}
\label{fig:costfunction_singlemodel}
\end{figure}

Fig.~\ref{fig:costfunction_singlemodel_topo} displays a heatmap of
this function for $n_{\theta} = n_{\psi} = 64$. We purposefully
construct this cost function to have a long ``ridge'' of local minima,
along with a relatively isolated global minimum. The reason for this
is to introduce topography that might be missed by a surrogate
constructed from too sparse of a data sampling. Given that
$\theta_{\text{true}}$ coincides with the region around the global
minimum, a relatively high degree of resolution in
$\mathcal{\widetilde{J}}$ is needed to accurately predict the optimal
strategy. We wish to explore how the isolated global minimum
necessitates our adaptive sampling and refinement strategy when using
a surrogate constructed from sparse data.

Regarding experimental design, we assume the following possible
experiment space:

\begin{equation}
  \begin{aligned}
    X &= \lbrace x_1 , \dots , x_{16} \rbrace \;\;\; \text{where} \; x_k = \theta_{4k} \\
    Y &= \lbrace 0 , 1 \rbrace \\
    \rho( y=1 | x , \theta ) &\sim  \; \text{exp}\left[ -\frac{1}{2} (x - \theta)^2 / \sigma_x^2 \right] \; , \; \sigma_x = \frac{1}{8}n_{\theta} \\
    \rho( y=0 | x , \theta)  &= 1 - \rho( y=1 | x , \theta )
  \end{aligned}
\end{equation}

One of the first issues that should be investigated is how our
adaptive sampling criterion tends to perform on this
problem. Fig.~\ref{fig:costfunction_singlemodel_sensitivity} shows the
sensitivity computed with our leave-one-out strategy, averaged over 64
Monte Carlo realizations. For comparison,
Fig.~\ref{fig:costfunction_singlemodel_DJ} shows the magnitude of the
gradient of the cost function, approximated by finite-differencing GP
surrogates and averaging over 64 Monte Carlo realizations. It is interesting
that these results seem to indicate that even locally flat regions
(e.g., very close to the global minimum) yield high sensitivity using
the leave-one-out strategy. We believe this is due to sparse sampling:
if, for instance, you only have one training point near the global
minimum, removing that training point will strongly affect the
surrogate, even though the local gradient is mild. For this reason, we
might speculate that the leave-one-out strategy might provide a better
measure of local ``importance'' than a gradient calculation when the
important topography of the cost function in question involves isolated
minima.

\begin{figure}[h!]
  \begin{minipage}{0.32\columnwidth}
    \centering 
    \subfloat[Initial surrogate $\widetilde{\mathcal{J}}$]{ \includegraphics[trim=0cm 0cm 0cm 0cm, clip, width=0.90\linewidth]{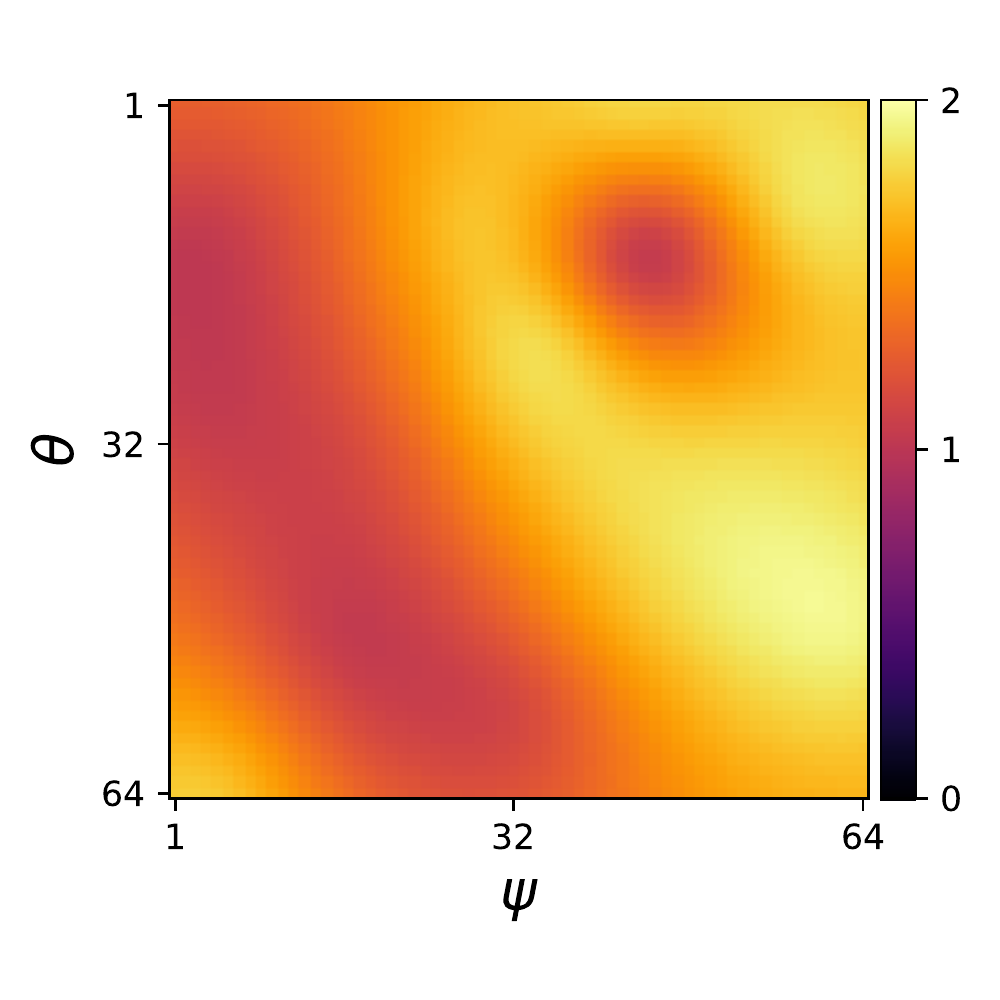} \label{fig:initialgpsurrogate_a} }\\
    \subfloat[Refined surrogate $\widetilde{\mathcal{J}}$]{ \includegraphics[trim=0cm 0cm 0cm 0cm, clip, width=0.90\linewidth]{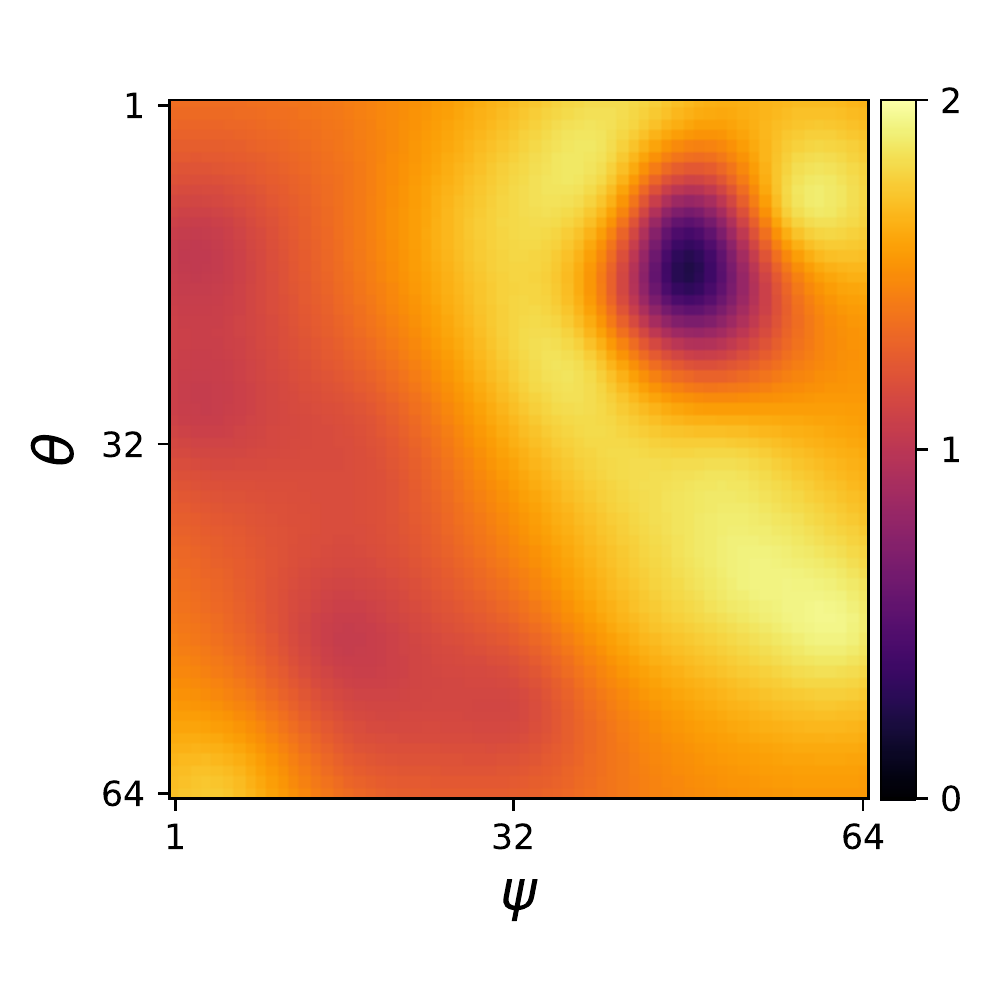} \label{fig:refinedgpsurrogate_a} }
  \end{minipage}
  \begin{minipage}{0.32\columnwidth}
    \subfloat[$ \left| \widetilde{\mathcal{J}} - \mathcal{J} \right|_1$]{ \includegraphics[trim=0cm 0.0cm 0cm 0cm, clip, width=0.90\linewidth]{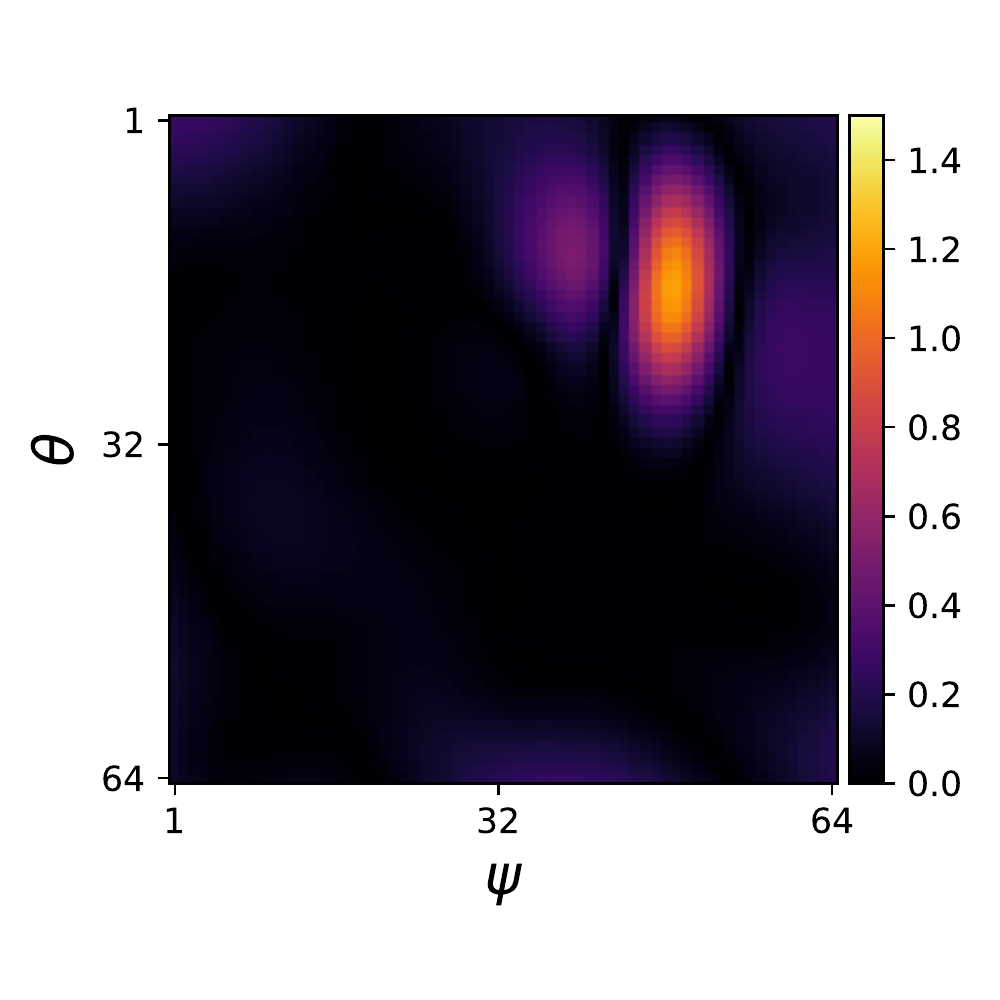} \label{fig:initialgpsurrogate_b} }\\
    \subfloat[$ \left| \widetilde{\mathcal{J}} - \mathcal{J} \right|_1$]{ \includegraphics[trim=0cm 0.0cm 0cm 0cm, clip, width=0.90\linewidth]{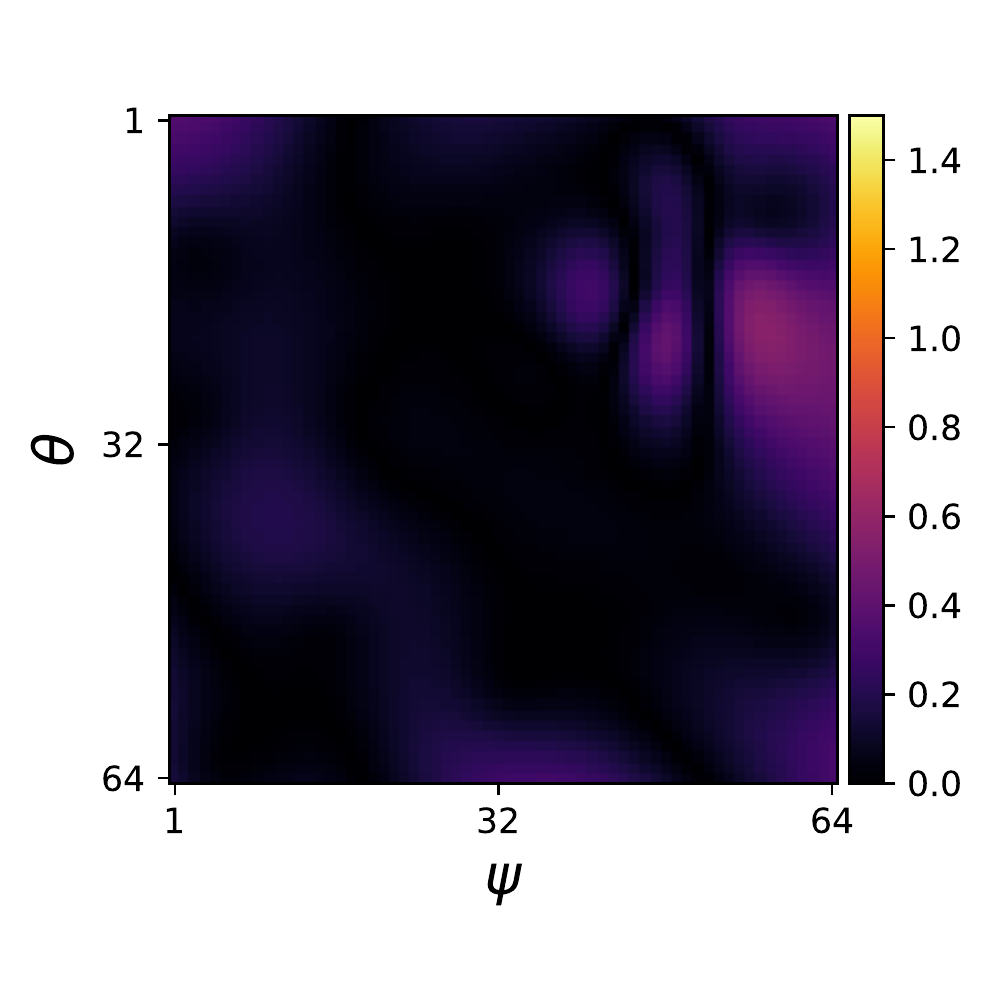} \label{fig:refinedgpsurrogate_b} }
  \end{minipage}
  \begin{minipage}{0.32\columnwidth}
    \subfloat[Training points]{ \includegraphics[trim=0cm 0cm 0cm 0cm, clip, width=0.80\linewidth]{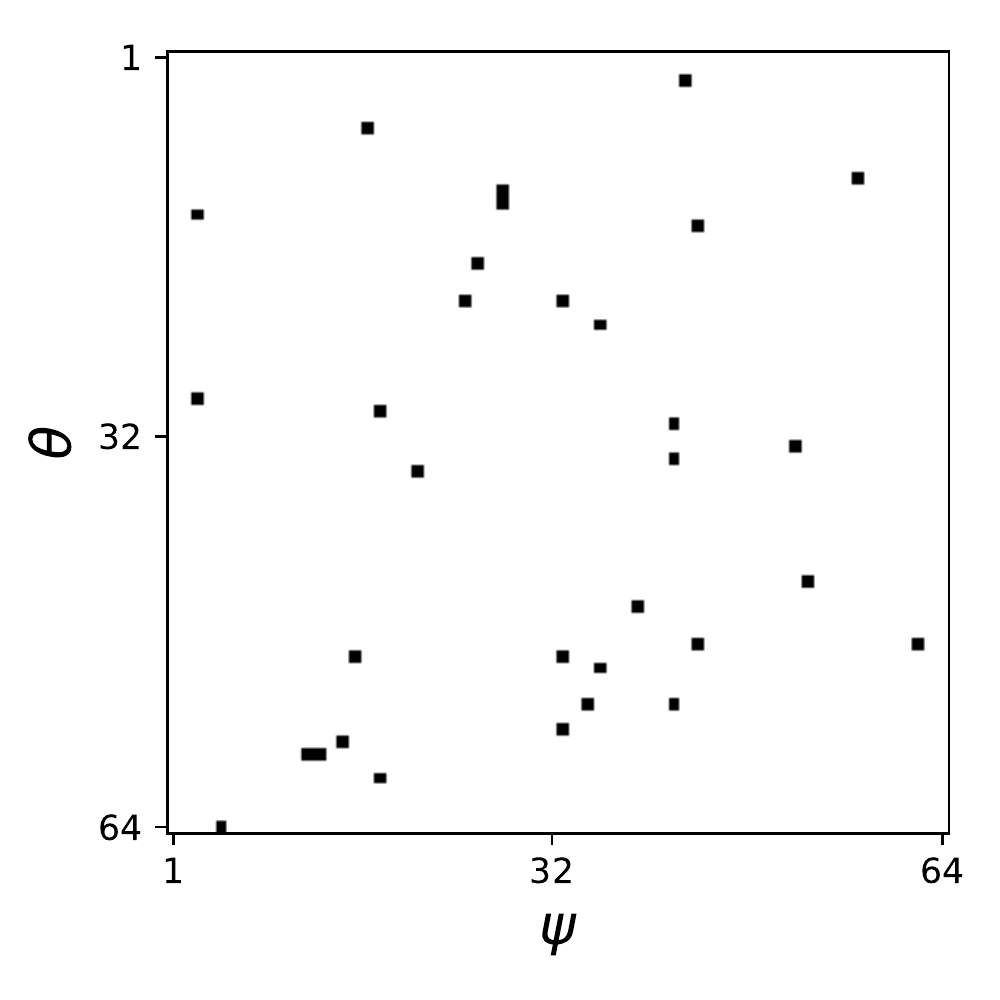} \label{fig:initialgpsurrogate_c} } \vspace*{0.9cm}
    \subfloat[New training points]{ \includegraphics[trim=0cm 0cm 0cm 0cm, clip, width=0.80\linewidth]{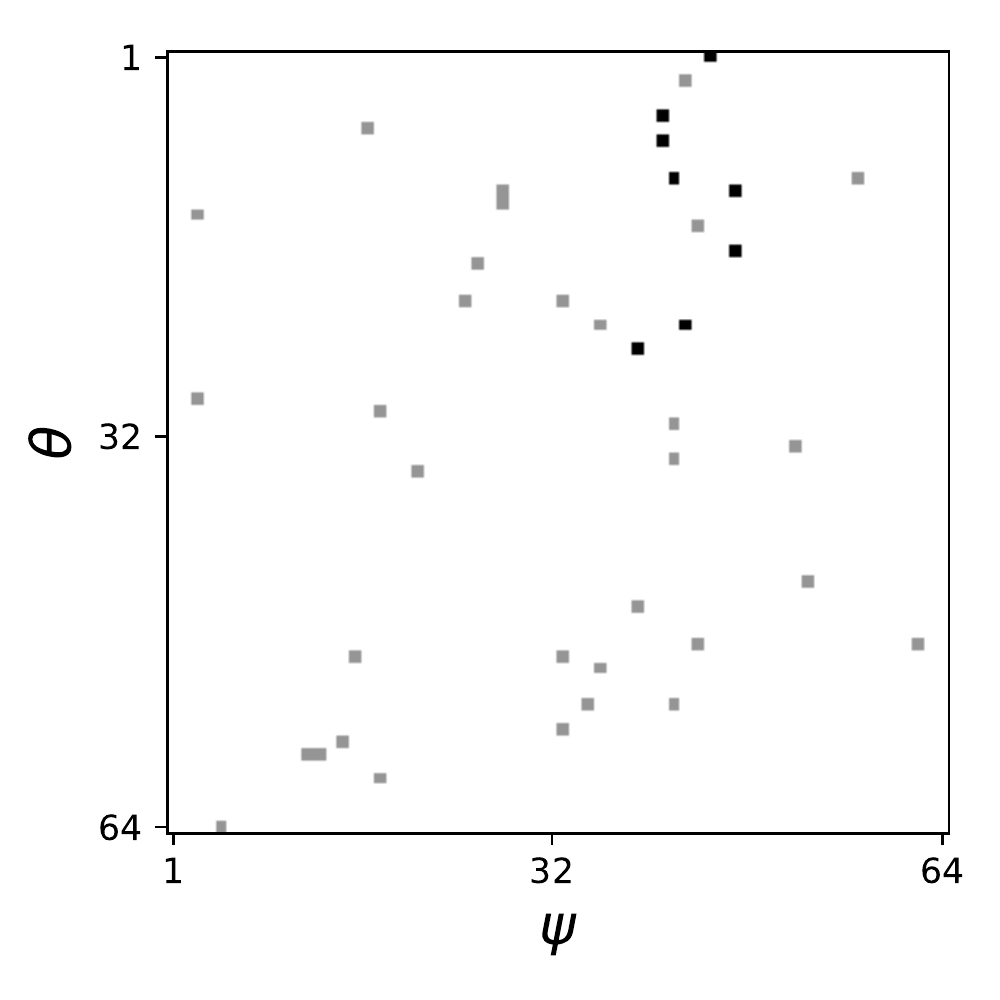} \label{fig:refinedgpsurrogate_c} }
  \end{minipage}
\caption{Adaptive-refinement procedure on an example GP surrogate. Top row: initial surrogate. Bottom row: Adaptively-refined surrogate. Original points are displayed in gray; new points are displayed in black.}
\label{fig:gpsurrogate}
\end{figure}

\begin{figure}[h!]
\centering
\includegraphics[width=0.90\textwidth]{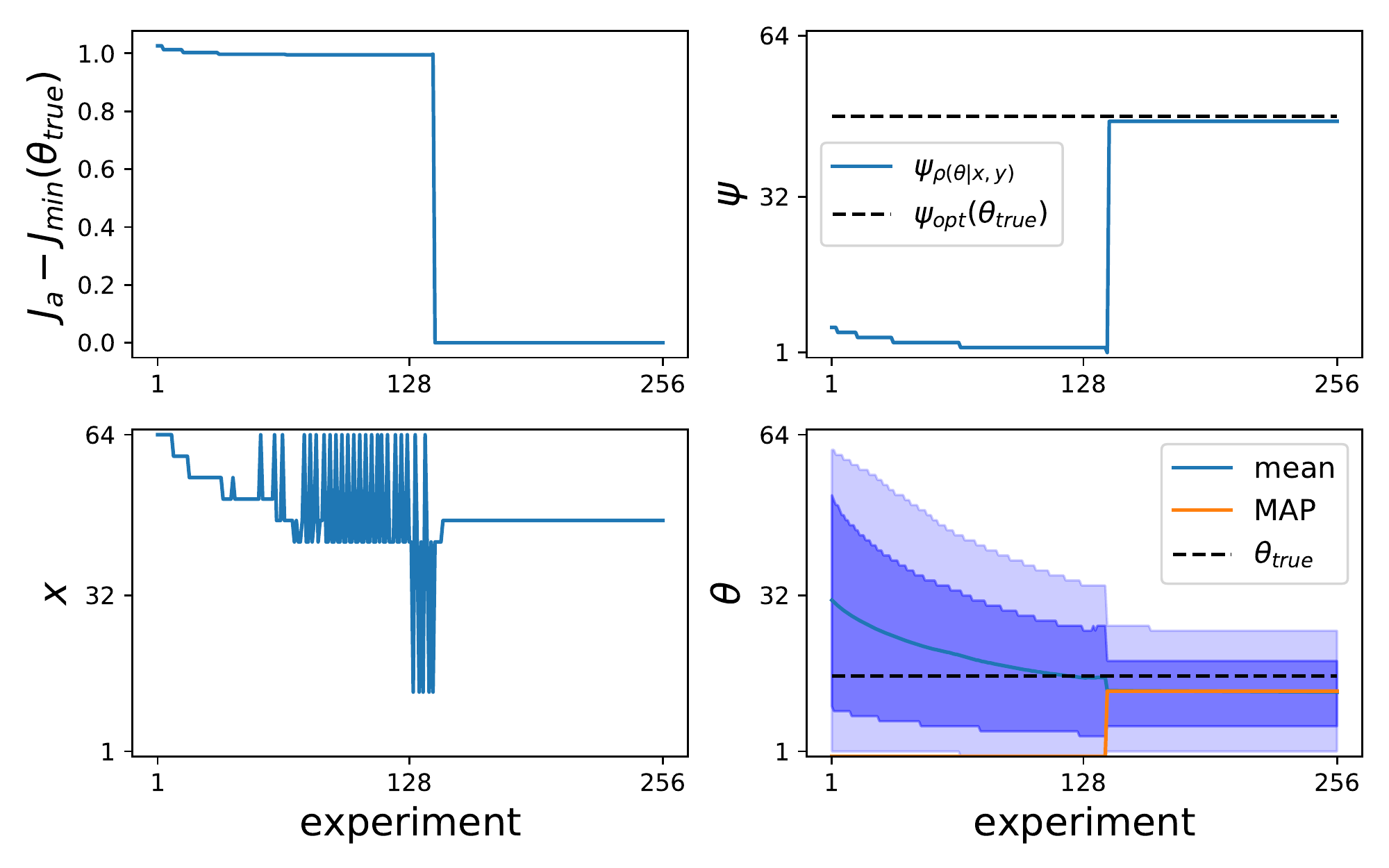}
\caption{Single cost model results: example MOCU results using adaptive sampling. Bottom-right: ``MAP'' stands for ``Maximum A-Posteriori'', and the shaded regions correspond to 68- and 95-\% confidence regions. $J_a$ denotes the adaptive MOCU cost.}
\label{fig:singlemodel_example_mocu_results}
\end{figure}

Fig.~\ref{fig:gpsurrogate} shows an example of a surrogate that has
been produced through adaptive refinement for this problem. We observe
good selective refinement near the global minimum as
desired. Fig.~\ref{fig:singlemodel_example_mocu_results} shows an
example of the OED results that are obtained using MOCU with adaptive
surrogate construction for this problem. Initially, there is high
entropy in the uncertainty class distribution, and the cost function
is approximated with a coarse surrogate, so the computed policy is far
from the optimal one. Eventually, once enough experiments have been
conducted, the distribution over the uncertainty class collapses close
to the correct value, and this triggers an adaptive refinement in the
regions near $\theta_{\text{true}}$. This results in training samples being
added near the global minimum, which in turn results in a much closer
approximation of the optimal policy.

Of course, the results shown in Fig.~\ref{fig:gpsurrogate} and
Fig.~\ref{fig:singlemodel_example_mocu_results} are just one
example. Because random sampling is involved in both constructing the
initial surrogate and refining it, we need to conduct multiple
independent Monte Carlo realizations to marginalize these effects
out. Further, we need to do this for both the adaptive and
non-adaptive schemes in order to confirm that the adaptive scheme
actually provides some benefit. To highlight the effects of
adaptivity, we also ensure that the total number of sample points used
in the adaptive case never exceeds the total number used in the
non-adaptive case. Thus, we use $48$ training points for each
non-adaptive run; in the adaptive case, we construct the initial
surrogate with $32$ training points and allow the option of adaptively
refining twice with 8 points each time (for a total no higher than
$48$).

\begin{figure}[h!]
\centering
\begin{multicols}{3}
    \includegraphics[trim=0cm 0cm 0cm 7.5cm, clip, width=0.95\linewidth]{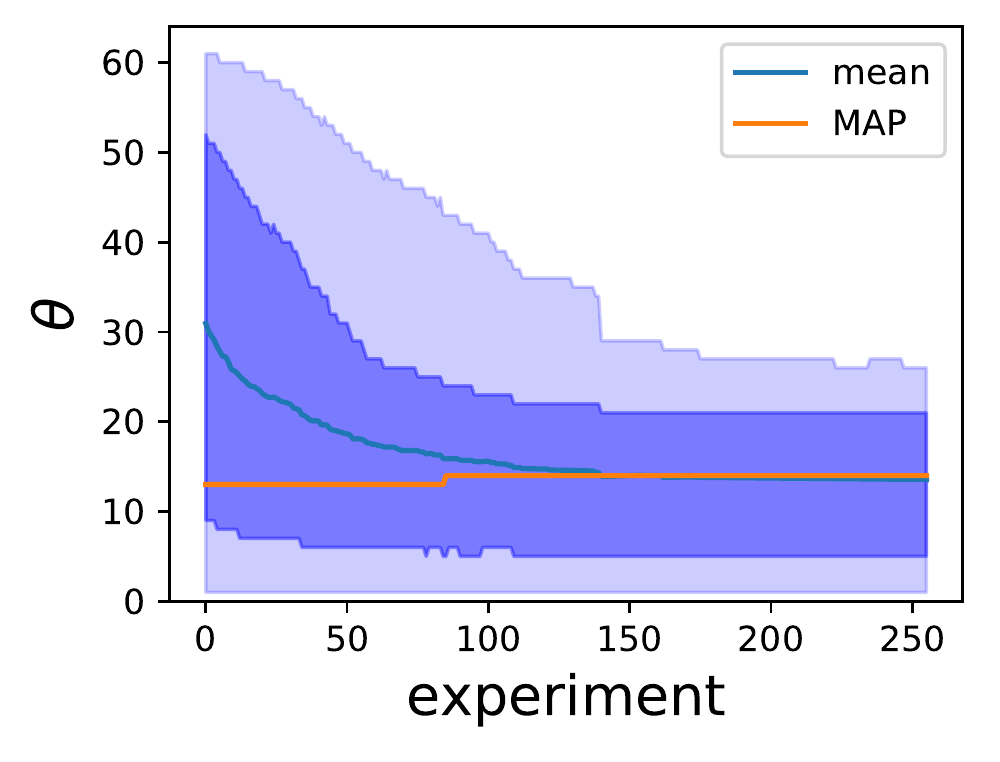}\par 
    \includegraphics[trim=0cm 0cm 0cm 7.5cm, clip, width=0.95\linewidth]{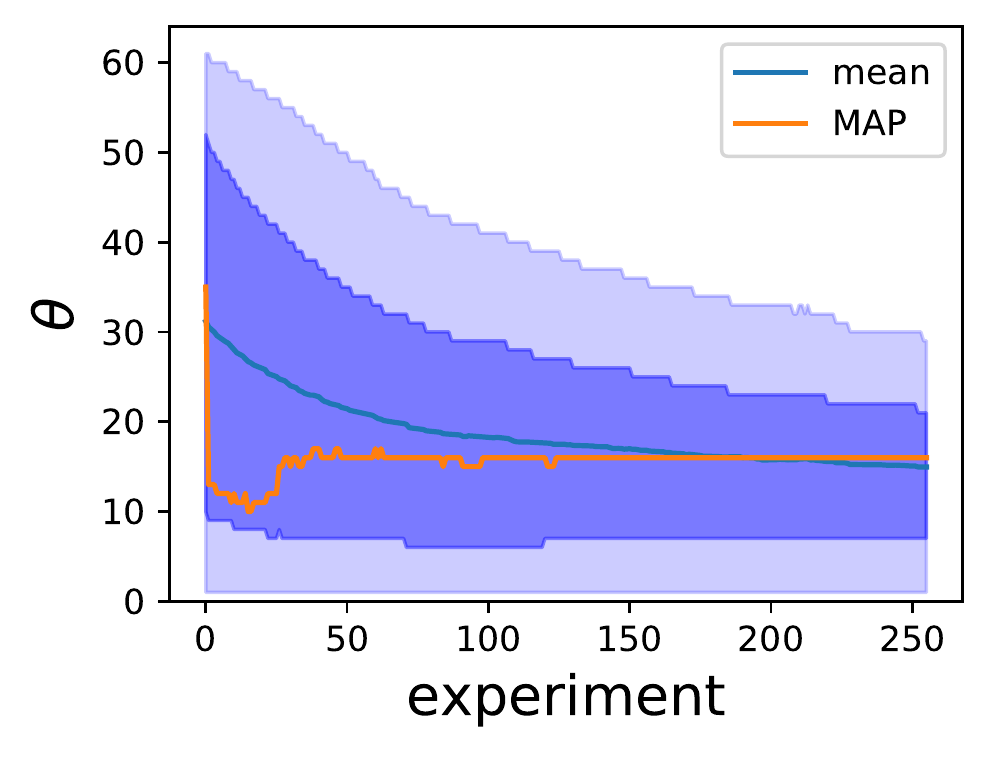}\par 
    \includegraphics[trim=0cm 0cm 0cm 7.5cm, clip, width=0.95\linewidth]{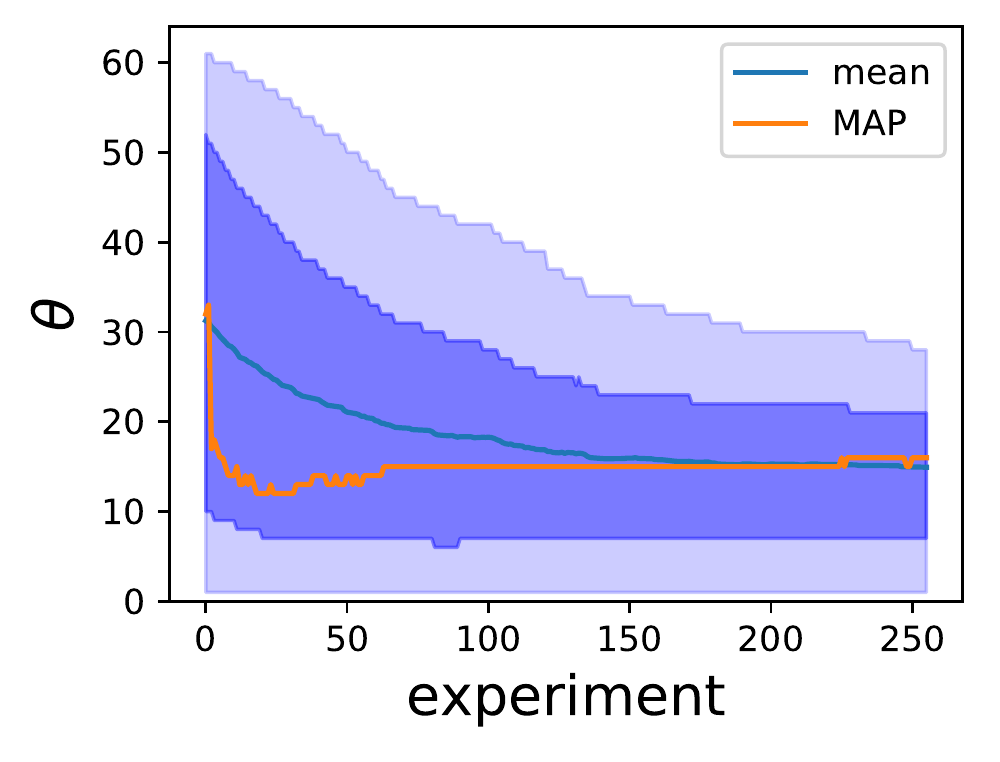}\par 
\end{multicols}
\caption{Single cost model Monte Carlo results: convergence of the uncertainty class distribution for standard MOCU (left) and non-adaptive (middle) and adaptive (right) sparse MOCU. ``MAP'' stands for ``Maximum A-Posteriori''.}
\label{fig:singlemodel_mc_rhotheta}
\end{figure}

Fig.~\ref{fig:singlemodel_mc_rhotheta} displays the evolution of
$\rho(\theta | x,y)$ with experiment for both the adaptive and non-adaptive
cases and for full MOCU (i.e., MOCU conducted with access to the full
cost function), averaged over 128 Monte Carlo realizations per method. We observe that,
on average, all three methods give roughly the same convergence
properties in $\rho(\theta | x,y)$. This suggests that in the MOCU method, epistemic
uncertainty in $\rho(\theta | x,y)$ is reduced at a rate that only weakly
depends on the accuracy of the computed cost function
$J(\theta,\psi)$. This may be a somewhat surprising observation, but
it is one of the foundational insights to our approach. If indeed we
will have to ``wait'' until a certain amount of experiments have been
conducted before we have a good estimate of $\theta_{\text{true}}$ --
and if this process is relatively unaffected by the use of an
approximate $\widetilde{\mathcal{J}}$ -- then it should be possible to
use a cheaply-constructed surrogate for the cost function, and
adaptively refine it only once we have conducted enough experiments.

\begin{figure}[h!]
    \centering 
    \includegraphics[trim=0cm 0cm 0cm 0cm, clip, width=0.6\linewidth]{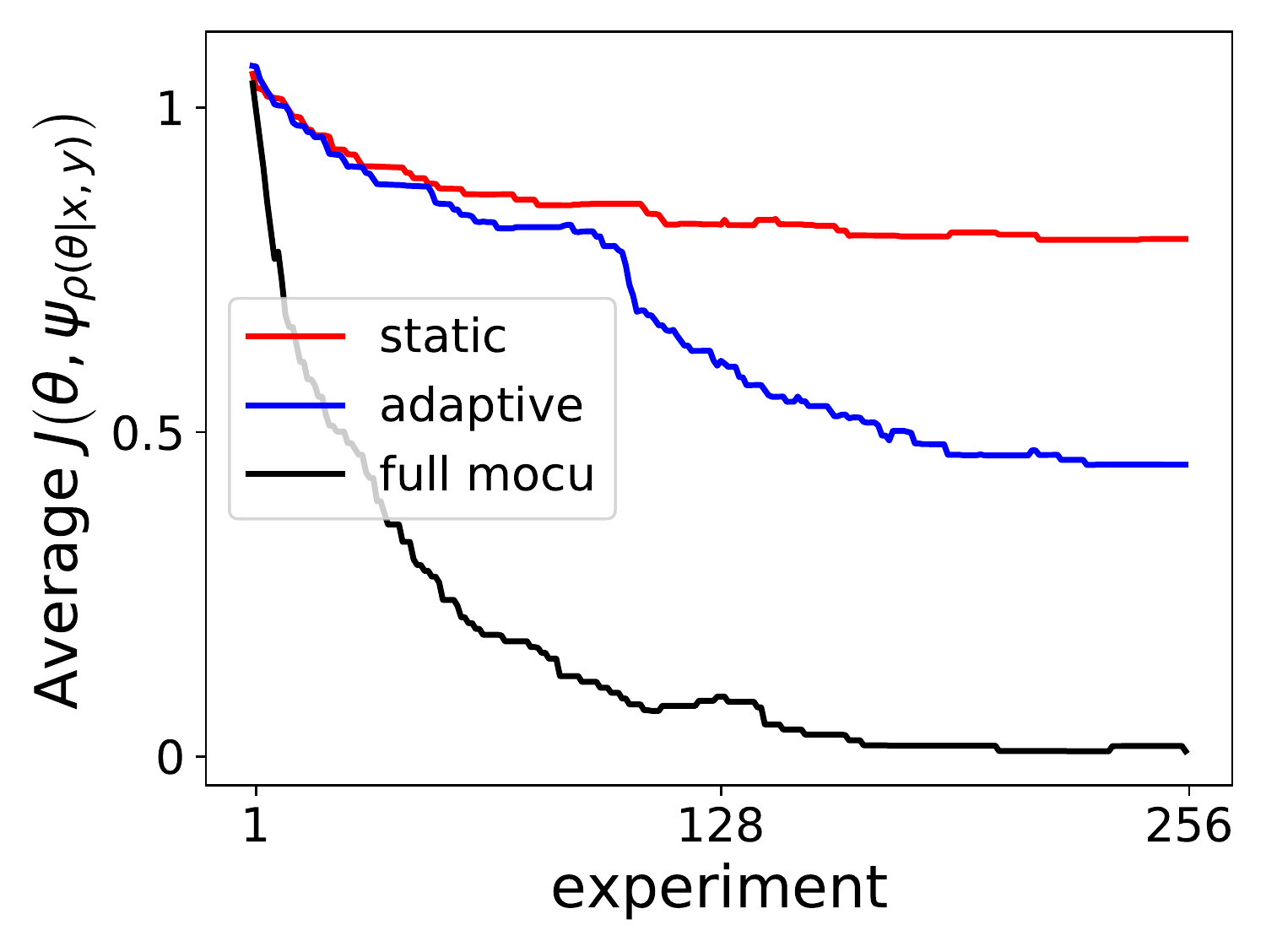} \label{fig:singlemodel_mc_results_a}
\caption{Single cost model Monte Carlo results. Averages at each experiment are computed over the 128 Monte Carlo samples.}
\label{fig:singlemodel_mc_results}
\end{figure}

Notwithstanding similarities in the evolution of $\rho(\theta | x,y)$
among all methods, we do observe significant differences in the
optimal policy calculation. Fig.~\ref{fig:singlemodel_mc_results}
displays results related to the optimal policy recommendation. In this
figure, we compute the average of the $J(\theta,
\psi_{\rho(\theta|x,y)})$ over the 128 Monte Carlo realizations at
each experiment. We see that full MOCU obviously has the highest
performance in locating the correct optimal policy (because it has
access to the full cost function $J(\theta,\psi)$). The rate at which
it converges roughly mirrors the convergence of $\rho(\theta |
x,y)$. This agrees with intuition in that the performance of full MOCU
is limited by the rate at which one can reduce epistemic uncertainty
through successive experiments. Regarding the two approximate schemes,
we see that the adaptive method tends to outperform the non-adaptive
one after a sufficient number of experiments have been performed (in
this case, about 100). This is because the distribution of the
uncertainty class has, on average, converged sufficiently by around
100 experiments that the adaptive refinement is triggered.

The comparison between full MOCU and the approximate schemes is
interesting from a practical perspective, and the results suggest a
trade-off between time spent evaluating the cost matrix and doing
experiments. On one hand, each of the adaptive realizations used between
$32$ and $48$ data points (i.e., the number of points in $\mathcal{P}$
used to construct the surrogate), whereas a full MOCU realization required
$n_{\theta} n_{\psi} = 64^2 = 4096$ evaluations of the cost
function. On the other hand, we see from
Fig.~\ref{fig:singlemodel_mc_results} that full MOCU requires
roughly an order of magnitude fewer experiments to achieve the same
level of accuracy as adaptive MOCU required after 256
experiments. Even though the evolution of the uncertainty class
distribution is almost identical for the two methods, full MOCU has
access to the entire cost matrix and thus is more robust against
larger amounts of uncertainty. The choice of whether to use our
adaptive, surrogate-aided method over full MOCU will depend on how
costly experiments are relative to evaluations of the cost matrix. If
the cost matrix is very expensive to evaluate (e.g., a multiscale
physics code) and experiments are inexpensive (e.g., there is a
pre-existing repository of collected historical data and
measurements), then one should consider using our approximate
methods. If instead experiments are prohibitively expensive (e.g., a
new exotic material must be fabricated and tested in a laboratory
setting to yield a single measurement), then one should consider using
regular, full MOCU.

\subsection{Multiple Model Cost Function}

Here, we extend the example from the previous section to a setting
where there are two cost functions: one that computes an
approximation of $J(\theta,\psi)$, and another that computes
$J(\theta,\psi)$ exactly. The motivation for introducing this is the
problem of multifidelity models~\cite{peherstorfer2018survey}. It is often the case that there are
multiple computer models that all predict the same
quantity-of-interest, but do so with computational expense inversely
proportional to accuracy. In such a setting, we might wonder whether
we could use the ``cheap'' approximate model to construct the initial
surrogate, and then use samples from the ``expensive'' accurate model
to refine the surrogate.

\begin{figure}[h!]
\centering
\begin{multicols}{2}
    \includegraphics[trim=4cm 0cm 5cm 0cm, clip, width=0.75\linewidth]{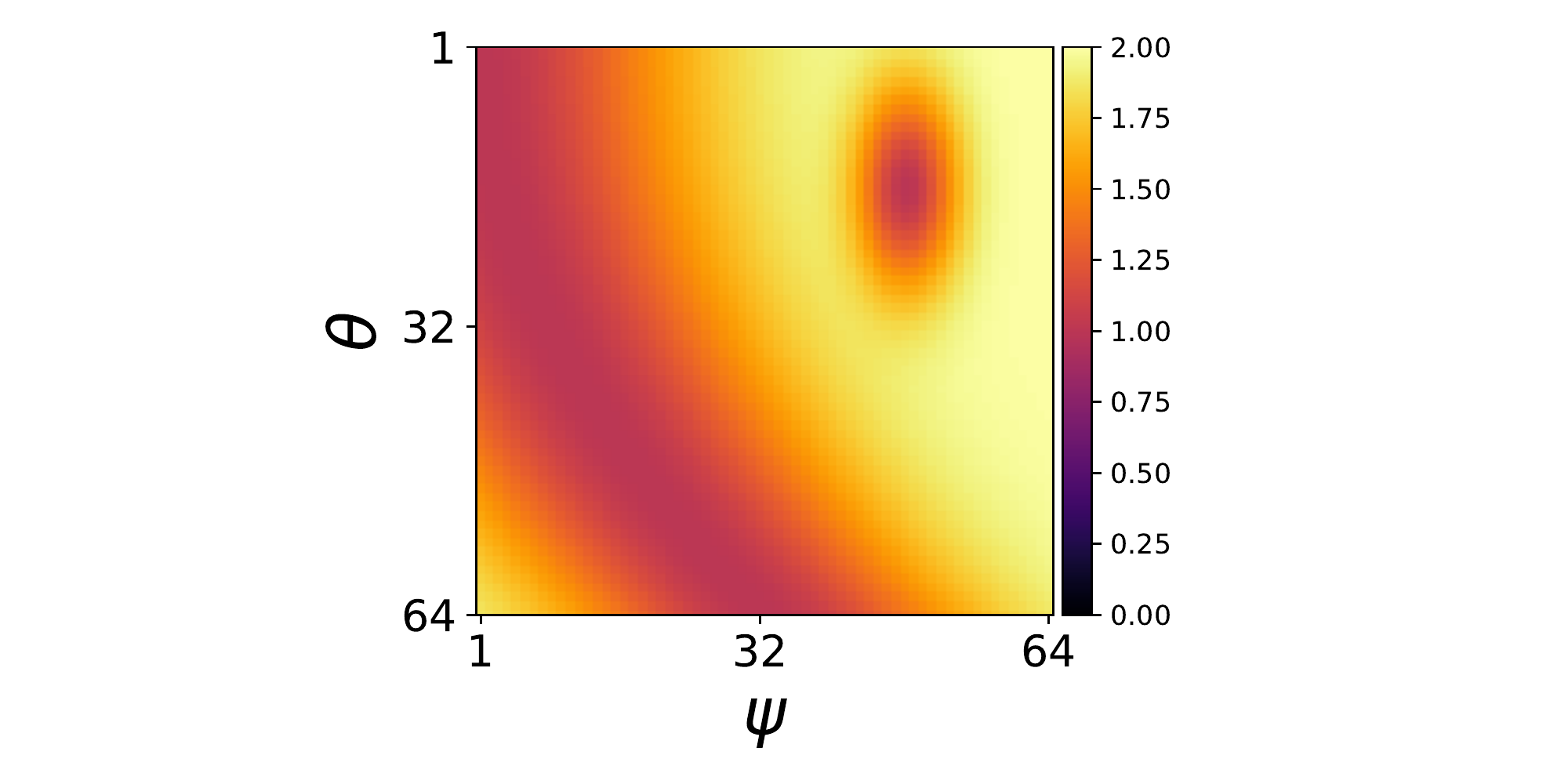}\par 
    \includegraphics[trim=4cm 0cm 5cm 0cm, clip, width=0.75\linewidth]{figures/multifidelity_fine_costfunction.pdf}\par 
\end{multicols}
\caption{Left: Coarse model cost function topography. Right: Fine model cost function topography.}
\label{fig:costfunction}
\end{figure}

Fig.~\ref{fig:costfunction} displays the cheap/inaccurate and
expensive/accurate cost functions used in this example. As can be
seen, it is impossible to distinguish between the isolated minimum and
the ridge minima on the basis of cost computed by the cheap model. It
is therefore impossible for MOCU to compute the correct policy based
on samples from the cheap model only; selective refinement from the
expensive model is needed in the vicinity of the global minimum.

\begin{figure}[h!]
\centering
    \includegraphics[trim=0cm 0cm 0cm 0cm, clip, width=0.6\linewidth]{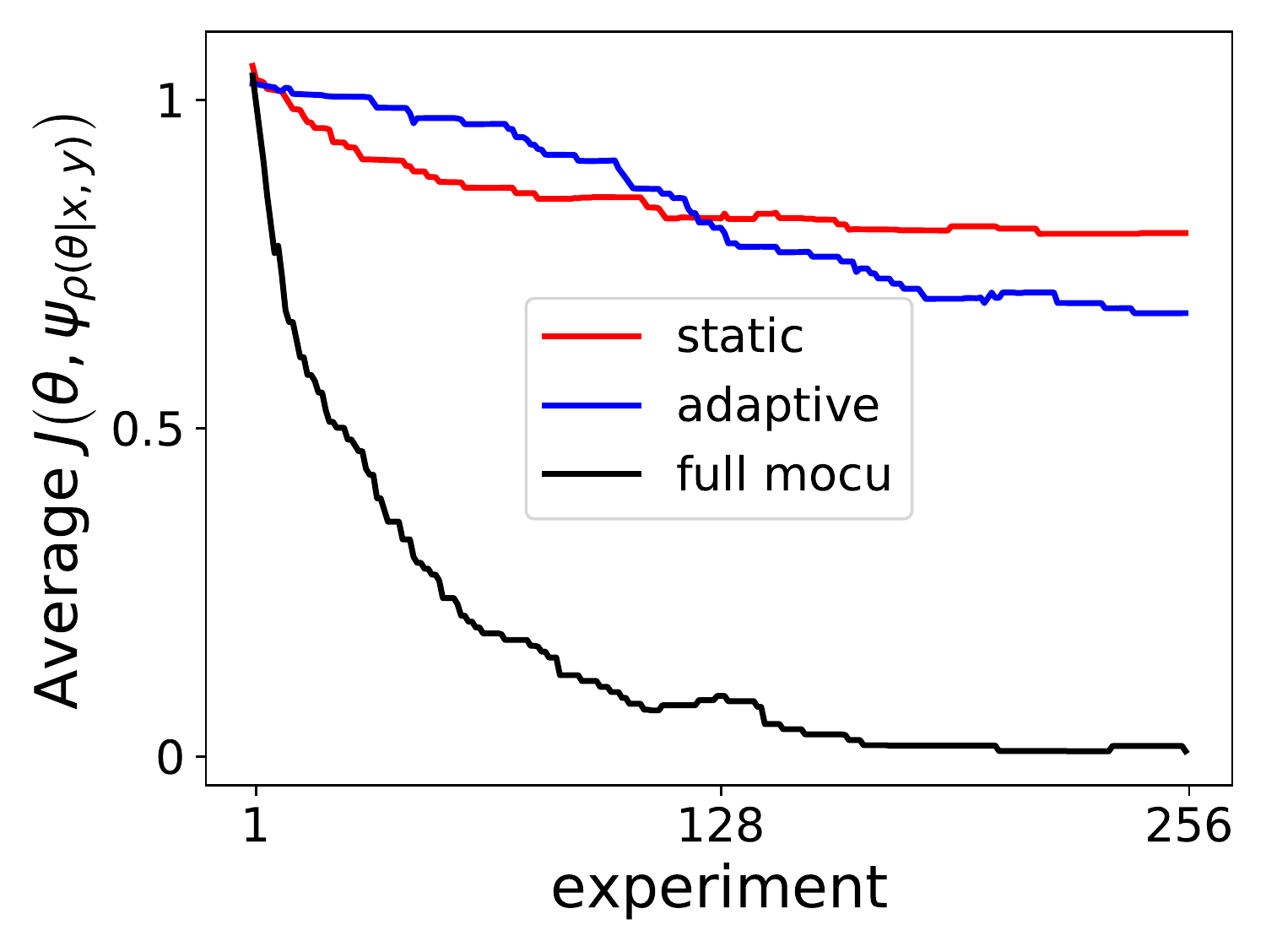}\par 
\caption{Multifidelity Monte Carlo results. Averages at each experiment are computed over the 128 Monte Carlo samples.}
\label{fig:multiplemodel_mc_results}
\end{figure}

Fig.~\ref{fig:multiplemodel_mc_results} reports Monte Carlo results
for this case, using both adaptive and non-adaptive sampling and full
MOCU. Note that both full MOCU and the non-adaptive surrogate-aided
methods are done with respect to the expensive model only, so the
results for those two cases are the same as in
Fig.~\ref{fig:singlemodel_mc_results}. Given the nature of this
problem, we are mostly interested in how the adaptive algorithm
compares with its non-adaptive counterpart. As has been noted, in the
non-adaptive case, all training samples were drawn from the expensive
model; in contrast, most samples for the adaptive case were drawn from
the cheap model, with only a minority selectively drawn from the
expensive model. And, as before, the total number of training samples
for the adaptive case never exceeds the number used for the
non-adaptive one (though it could be less). Even with these
disadvantages, we observe that there is still
statistically significant improvement in the optimal policy
calculation using the adaptive method, which begins to become apparent
after about 128 experiments.

\subsection{Application: Coupled Spring-Mass-Damper Design}

Here, we take a first step towards applying our methods to the design
and control of physical systems. The problem setting will be a
controls problem for a linear system. Specifically, we consider the
following coupled spring-mass-damper system:

\begin{equation}
  \begin{aligned}
    m_1 \ddot{x}_1 &= -k_1  x_1 - k_2 ( x_1 - x_2 ) - \delta_1  v_1 - \delta_2 ( v_1 - v_2 ) \\
    m_i \ddot{x}_i &= -k_i( x_i - x_{i-1} ) - k_{i+1}( x_i - x_{i+1} ) \\ & \;\;\;\; - \delta_i( v_i - v_{i-1} ) - \delta_{i+1}( v_i - v_{i+1} ) \;\;\; , \;\;\; i=2 \dots n-1\\
    m_n \ddot{x}_n &= -k_n( x_n - x_{n-1} ) - \delta_n( v_n - v_{n-1} )
  \end{aligned}
\end{equation}

We can recast this in first order form, along with an output QOI and control signal:

\begin{equation}
  \begin{aligned}
  \frac{d}{dt}
  \begin{bmatrix}
    \mathbf{x} \\
    \mathbf{v}
  \end{bmatrix} &=
  \begin{bmatrix}
    0 & \mathds{1} \\
    A_2 & A_3
  \end{bmatrix}
  \begin{bmatrix}
    \mathbf{x} \\
    \mathbf{v}
  \end{bmatrix} +
  \begin{bmatrix}
    \mathbf{0} \\
    \mathbf{b}
  \end{bmatrix}
  u \\
  y &= \begin{bmatrix}
    \mathbf{c}^T & \mathbf{0}^T
  \end{bmatrix}
  \begin{bmatrix}
    \mathbf{x} \\
    \mathbf{v}
  \end{bmatrix} \;\;\; .
  \end{aligned}
  \label{eqn:linsys_springs}
\end{equation}

In Eqn.~\ref{eqn:linsys_springs}, $\mathbf{x} = [ x_1 , \dots , x_n ]$
denote the displacements from equilibrium of the $n$-springs, and
$\mathbf{v} = [ v_1 , \dots , v_n ]$ denotes the 1-D velocities of
these springs. In this problem, we let $\mathbf{b} = [1 , 0 , \dots ,
  0]$ and $\mathbf{c} = \left[ \frac{1}{n} , \dots , \frac{1}{n} \right]$,
i.e. our input signal is a force that affects the first mass in the
chain, and we are measuring the average spring displacement across
the entire chain.

Our goal is to find the input sinusoidal signal that maximizes our
output amplitude. This occurs at a resonant frequency of our
system. For fixed system parameters, we can compute this from
examination of the Bode magnitude plot of the transfer function from
$U(s)$ to $Y(s)$, where $U(s)$ and $Y(s)$ are the Laplace transforms
of $u(t)$ and $y(t)$. It is well-known that this transfer function --
which we will denote as $H(s) = \frac{Y(s)}{U(s)}$ -- can be computed as
$C(s \mathds{1} - A)^{-1} B$, where $(A,B,C)$ are the state-space
matrices in Eqn.~\ref{eqn:linsys_springs}.

In our problem, we allow for uncertainty in the $n$ spring
coefficients, and we wish to compute that sinusoidal forcing frequency
that best maximizes the output amplitude, on average across this
uncertainty. We note that one could use methods from optimal/robust
control to solve this problem, instead of MOCU. Our goal in
introducing this problem is simply to demonstrate our methods on a
minimally-complex physical system.

\begin{figure}[h!]
  \begin{minipage}{0.48\columnwidth}
    \centering 
    \subfloat[Bode magnitude]{ \includegraphics[trim=0cm 0cm 0cm 0cm, clip, width=0.92\linewidth]{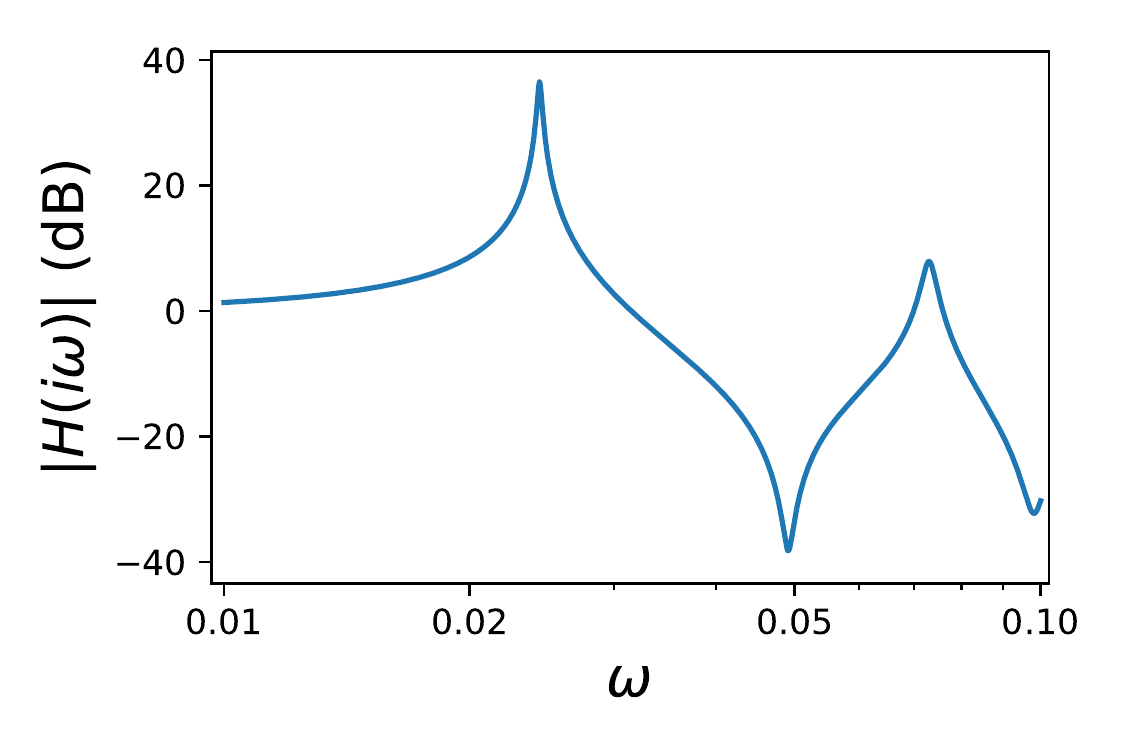} \label{fig:bode_and_J_springs_a} }
  \end{minipage}
  \begin{minipage}{0.48\columnwidth}
    \subfloat[Cost function]{ \includegraphics[trim=4cm 0cm 5cm 0cm, clip, width=0.75\linewidth]{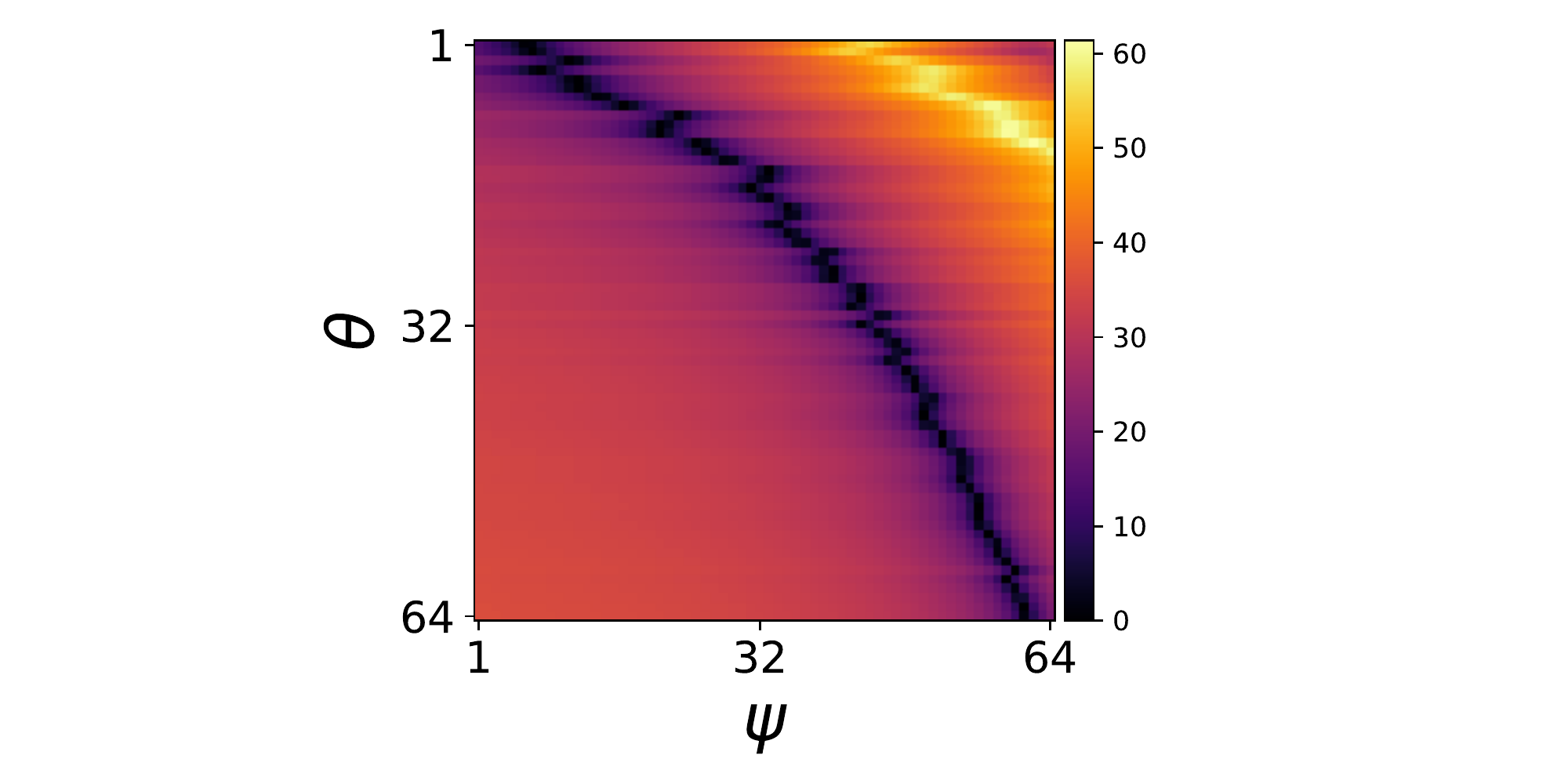} \label{fig:bode_and_J_springs_b} }
  \end{minipage}
\caption{Left: Bode magnitude plot for coupled spring-mass-damper system with spring parameters $k_i = 1 \; \forall i=1 \dots n$. Right: Example cost function for the spring-mass-damper system.}
\label{fig:bode_and_J_springs}
\end{figure}

Fig.~\ref{fig:bode_and_J_springs_a} plots the Bode magnitude plot for
the system with $k_i = 1 \;\forall i=1 \dots n$, in the range $\omega
\in [ 10^{-2} , 10^{-1} ]$. We can clearly see the presence of a
dominant resonant frequency as well as a secondary resonant frequency
in this range. In this problem, we consider $n=16$ springs. We set
$m_i = 1$ and $\delta_i = \frac{1}{8} \;\forall i=1 \dots n$. We set
the spring frequencies to $k_i = \theta_i = 0.1 + 0.9
\frac{i}{n_{\theta}} + \eta_i \;\forall i=1 \dots n$, where $\eta_i$
is drawn from a uniform distribution: $\eta_i \sim
\mathcal{U}[-0.1,0.1]$. Thus, our uncertainty class corresponds to
increasingly stiff systems (up to the white noise produced by
$\eta$). We also constrain the minimum allowable value of $k_i$ to
0.1. We take $n_{\theta} = 64$. The action set $\Psi = \lbrace 1 ,
\dots , n \rbrace$ sets the forcing frequency of the input sinusoid:
$\omega = \omega(\psi)$, where $\omega(\psi)$ is a log-space mapping
to the range $[0.03 , 0.1]$. Fig.~\ref{fig:bode_and_J_springs_b}
displays an example ground-true cost function $J(\theta,\psi)$ for
this problem. Further, we set the true value of the uncertain
parameter to be $\theta_{\text{true}} = \frac{3}{4}n$. We compute our
cost function as $J(\theta,\psi) = \text{max}_{\Psi}[ | H(\theta , i
  \omega(\psi) ) | ] - |H(\theta, i \omega(\psi))|$; that is, we
compute the cost for a specific $(\theta,\psi)$ pair as the deviation
in $|H|$ from the $\theta$-specific maximum in $|H|$. We note that
this is a move one cannot make in a truly adaptive application
(because it assumes one can calculate the maximum over $\psi$ for each
$\theta$), but we simply do it here for convenience. In practice, one
could circumvent this issue by instead considering an arbitrary
reference point, $J(\theta,\psi) = c - |H(\theta, i \omega(\psi))|$
for some $c \geq |H(\theta, i \omega(\psi))| \; \forall (\theta,
\psi)$, though that choice should not substantially affect the results
for this problem.

\begin{figure}[h!]
\centering
    \includegraphics[trim=0cm 0cm 0cm 0cm, clip, width=0.6\linewidth]{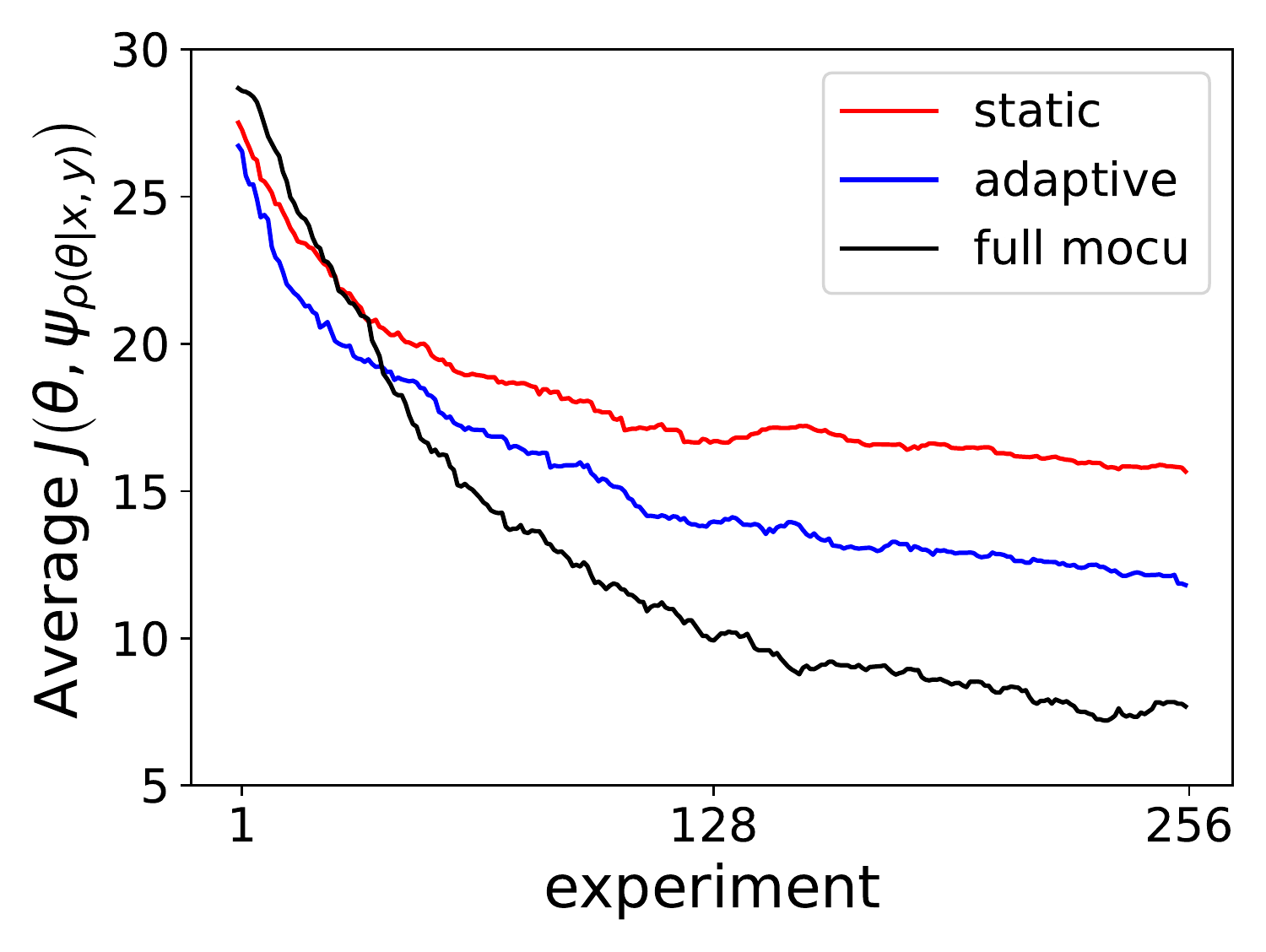}\par 
\caption{Coupled spring-mass-damper system Monte Carlo
  results. Averages at each experiment are computed over the 128 Monte
  Carlo samples.}
\label{fig:springs_mc_results}
\end{figure}

\begin{figure}[h!]
\centering
\includegraphics[width=0.90\textwidth]{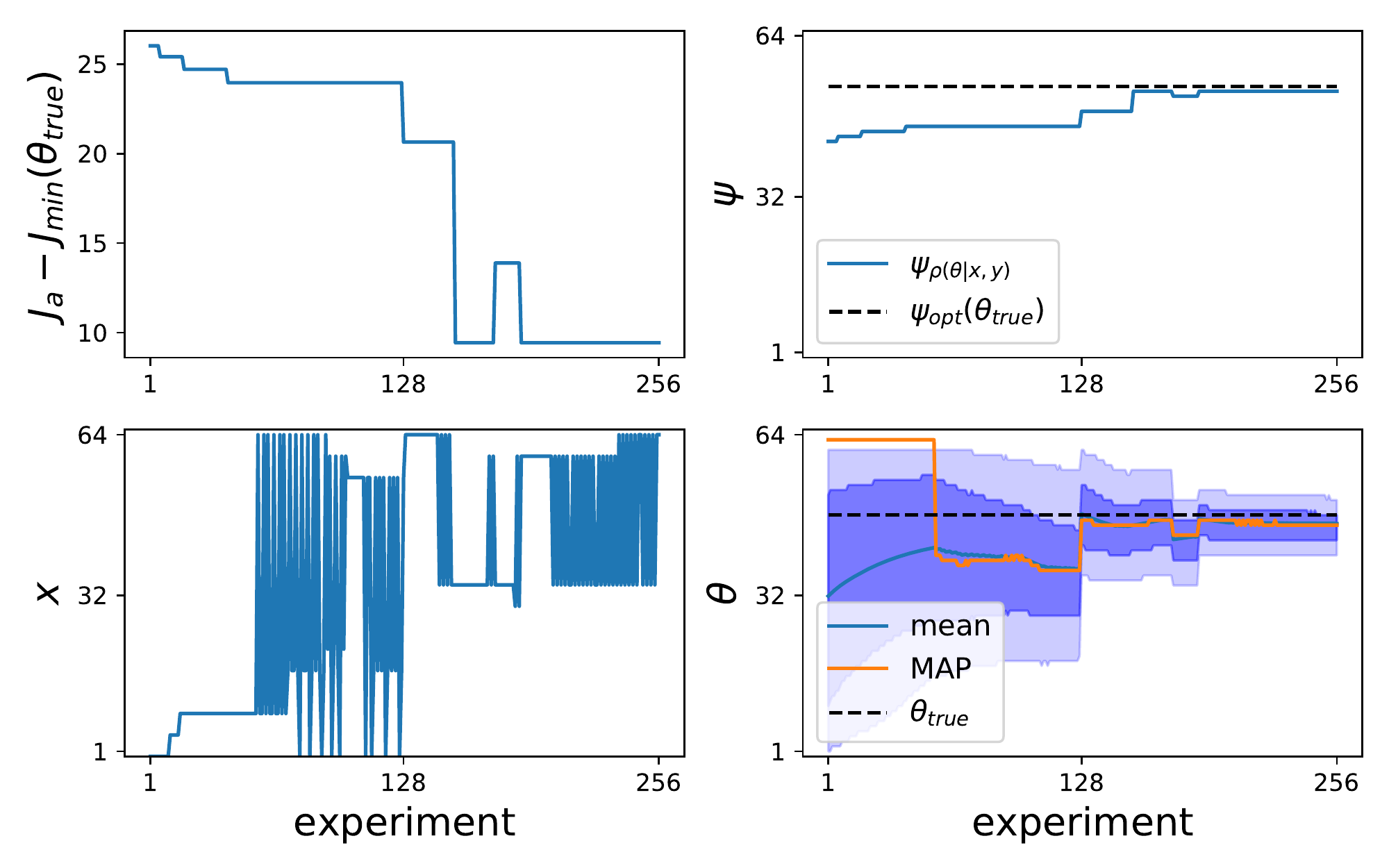}
\caption{Coupled spring-mass-damper system results: example MOCU results using adaptive sampling for the coupled spring-mass-damper system. ``MAP'' stands for ``Maximum A-Posteriori''. $J_a$ denotes the adaptive MOCU cost.}
\label{fig:springs_example_mocu_results}
\end{figure}

We apply the same three MOCU methods to this problem as in the
previous examples (i.e., non-adaptive and adaptive MOCU, with a
surrogate for the cost function, and full MOCU). Note that in this
problem, all surrogates (adaptive and non-adaptive) are initially
constructed with 32 training points, and the adaptive surrogates are
allowed two possible refinements of 8 points (for a total of 32 to 48
training points). Fig.~\ref{fig:springs_example_mocu_results} displays
an example run of adaptively refined, surrogate-approximated MOCU for
this problem, and Fig.~\ref{fig:springs_mc_results} displays the
statistical comparisons of the performance of full and approximate
MOCU schemes. In comparing the two approximate methods, we find a
predictable ranking: the adaptive method more accurately estimates the
optimal policy with higher probability than the non-adaptive
method. Regarding the full MOCU results, we observe the same trade-off
that we noted previously: full MOCU can attain the same level of
statistical accuracy as surrogate-approximated MOCU (or better), but
with fewer needed experiments. As before, the decision of whether to
use full MOCU or surrogate-approximated MOCU will depend on the
relative cost of computation versus experimentation.

\section{Conclusions}

The goal of this research was to propose a new strategy for
approximate OED for resource-constrained problems. Our motivation was
to make OED (via MOCU) tractable for settings where the computational
load needed to evaluate the design cost over the full set $\Theta
\times \Psi$ is prohibitively large. This could occur either because a
single evaluation of $J(\theta_i,\psi_j)$ is expensive for all $i,j$
pairs, or because the joint set $\Theta \times \Psi$ is of large size,
or both. Thus, our focus was on reducing the computationally-intensive
stages of MOCU. To do this, we introduced the idea of using a
surrogate model to approximate the cost function
$J(\theta,\psi)$. This surrogate is built initially from
sparsely-sampled data pairs drawn from $\Theta \times \Psi$, and it is
refined adaptively as more information is gathered about the
uncertainty class via the data-conditioned posterior $\rho(\theta |
x,y)$. We applied this method to several example problems and examined
its performance relative to a static, unrefined surrogate and full
MOCU. We conclude that the adaptive refinement generally improves the
performance of surrogate-driven approximate MOCU, but that the
decision on whether to use full MOCU versus our approximate methods
depends on the relative expense of evaluating the design cost versus
doing an experiment, on the feasibility of doing full MOCU, and on the
desired level of accuracy in optimal policy recommendations.

There are several avenues of further research that should be
investigated in the future. One is practical: the approximation
methods we have discussed should be applied to complex, ``real-world''
design problems. Another is theoretical: our methods are most useful
in the case that experiments are inexpensive relative to
compute-time. It would be useful to develop an extension that accounts
for the relative expense of experiments and computations, and suggests
a strategy that is weighted to account for this. For example, perhaps
the number of initial training points and the number of adaptive
refinements could be selected according to this criterion. One could
also investigate how performance varies with different methods for
adaptive selection (other than leave-one-out) and how the
hyperparameters for our scheme might be tuned.

\section*{Funding Sources}

This work was supported by the U.S. Department of Energy, Office of
Science, Office of Advanced Scientific Computing Research,
Mathematical Multifaceted Integrated Capability Centers program under
Award Q1 DE-SC0019393.

\section*{Acknowledgments}

We wish to thank Edward Dougherty, Byung-Jun Yoon and Nathan Urban for
their helpful and insightful discussions on MOCU.

\newpage
\bibliographystyle{siam_latex_template/siamplain}
\bibliography{paper_sparse_mocu}

\end{document}